\documentclass[10pt,a4paper,reqno]{amsart}
\usepackage{latexsym,amssymb,amsmath,amsthm,amsfonts,enumerate,verbatim,xspace,exscale}

\input xy
\xyoption{all}
\CompileMatrices
\UseComputerModernTips

\theoremstyle{plain}
\newtheorem{theorem}{Theorem}[section]
\newtheorem{lemma}[theorem]{Lemma}
\newtheorem{proposition}[theorem]{Proposition}
\newtheorem{corollary}[theorem]{Corollary}
\newtheorem*{Theorem}{Theorem}

\theoremstyle{definition}

\newtheorem{definition}[theorem]{Definition}

\newtheorem*{name}{\theoremname}
\newcommand{\theoremname}{testing}

\newtheorem{stepp}{Remark}[section]

\newtheorem{exx}[theorem]{Example}

\newenvironment{ex}
    {\begin{exx}}
    {\x\end{exx}}

\allowdisplaybreaks

\newcommand{\caps}[1]{\textup{\textsc{#1}}}

\providecommand{\bysame}{\makebox[3em]{\hrulefill}\thinspace}

\newcommand{\iso}{isomorphism\xspace}

\newcommand{\diffeo}{diffeomorphism\xspace}

\newcommand{\bb}{\mathbb}
\newcommand{\ov}[1]{\mbox{$\overline{#1}$}}
\newcommand{\up}{\upshape}

\newcommand{\x}{$\hfill\Box$}

\newcommand{\longto}{\longrightarrow}
\newcommand{\hookto}{\hookrightarrow}
\newcommand{\toto}{\twoheadrightarrow}

\def\vv<#1>{\langle#1\rangle}

\providecommand{\sign}{\mbox{$\text{\up{sign}}\,$}}

\newcommand{\dd}[2]{\mbox{$\frac{\partial #2}{\partial #1}$}}
\newcommand{\ddo}{\mbox{$\dd{t}{}|_{0}$}}

\newcommand{\om}{\omega}
\newcommand{\Om}{\Omega}
\newcommand{\var}{\varphi}

\newcommand{\lam}{\lambda}
\newcommand{\Lam}{\Lambda}

\newcommand{\wt}[1]{\mbox{$\widetilde{#1}$}}
\newcommand{\dwt}[1]{\mbox{$\widetilde{\widetilde{#1}}$}}

\newcommand{\poly}{\mbox{$\text{\up{Poly}}$}}

\newcommand{\bsc}{\mbox{$\bigsqcup$}}

\newcommand{\N}{\mbox{$\bb{N}$}}

\newcommand{\R}{\mbox{$\bb{R}$}}
\newcommand{\C}{\mbox{$\bb{C}$}}

\newcommand{\by}[2]{\mbox{$\frac{#1}{#2}$}}
\newcommand{\cinf}{\mbox{$C^{\infty}$}}
\providecommand{\set}[1]{\mbox{$\{#1\}$}}
\newcommand{\subeq}{\subseteq}

\newcommand{\fl}{\mbox{$\text{\up{Fl}}$}}

\newcommand{\X}{\mathfrak{X}}
\newcommand{\ld}{\mathcal{L}}
\newcommand{\curv}{\mbox{$\textup{Curv}$}}
\newcommand{\vf}{vector field\xspace}
\newcommand{\mf}{manifold\xspace}

\newcommand{\ver}{\mbox{$\textup{Ver}$}}
\newcommand{\hor}{\mbox{$\textup{Hor}$}}
\newcommand{\nor}{\mbox{$\textup{Nor}$}}

\newcommand{\ann}{\mbox{$\textup{Ann}\,$}}
\newcommand{\fix}[1]{\mbox{$\textup{Fix}(#1)$}}

\newcommand{\bra}[1]{\mbox{$\{#1\}$}}
\newcommand{\ham}[1]{\mbox{$\smash{\nabla^{\om}_{#1}}$}}

\newcommand{\spr}[1]{\mbox{$/\negmedspace/_{#1}$}}
\newcommand{\sporb}{\mbox{$/\negmedspace/_{\mathcal{O}}$}}
\newcommand{\momap}{momentum map\xspace}

\newcommand{\omkks}{\mbox{$\Om^{\mathcal{O}}$}}

\newcommand{\gu}{\mathfrak{g}}
\newcommand{\ho}{\mathfrak{h}}

\newcommand{\ko}{\mathfrak{k}}
\newcommand{\lo}{\mathfrak{l}}

\newcommand{\Ad}{\mbox{$\text{\upshape{Ad}}$}}
\newcommand{\ad}{\mbox{$\text{\upshape{ad}}$}}

\newcommand{\orb}{\mbox{$\mathcal{O}$}}

\newcommand{\WW}{\mbox{$\mathcal{W}$}}

\newcommand{\lorb}[1]{\mbox{${#1}_{(L)}^{\mathcal{O}}$}}
\newcommand{\dep}{\mbox{$\textup{depth}\,$}}
\newcommand{\bscorb}{\mbox{$\bsc_{q\in Q}\mathcal{O}\cap\ann\ko_{q}$}}
\newcommand{\ine}{\mbox{$\mathbb{I}$}}

\newcommand{\aklm}{Alekseevsky, Kriegl, Losik, Michor\xspace}

\title[Singular Poisson Reduction of Cotangent Bundles]
      {Singular Poisson Reduction of Cotangent Bundles}
\author{Simon Hochgerner\\
  and\\
  Armin Rainer}
\address{S.\ Hochgerner, Institut f\"{u}r Mathematik\\
  Universit\"{a}t Wien\\
  Nordbergstrasse~15\\
  A-1090 Vienna, Austria}
\email{simon.hochgerner@univie.ac.at} 
\urladdr{http://www.mat.univie.ac.at/$\thicksim$simon}
\address{A.\ Rainer, Institut f\"{u}r Mathematik\\
  Universit\"{a}t Wien\\
  Nordbergstrasse~15\\
  A-1090 Vienna, Austria}
\email{armin.rainer@univie.ac.at}
\urladdr{http://www.mat.univie.ac.at/$\thicksim$armin}
\thanks{This work is supported by Fonds zur F\"{o}rderung der
       wissenschaftlichen Forschung, Projekt P 17108-N04}
\date{August 24, 2005}

\begin{document}

\begin{abstract}
We consider the Poisson reduced space $(T^*Q)/K$
with respect to a cotangent lifted action.  
It is assumed that $K$ 
is a compact Lie group which acts by isometries on the Riemannian 
manifold $Q$ and that the action on $Q$ 
is of single isotropy type. 
Realizing $(T^*Q)/K$ as a
Weinstein space we determine the induced Poisson structure and  
its symplectic leaves. We thus extend the Weinstein construction for
principal fiber bundles to the case of surjective Riemannian
submersions $Q\toto Q/K$.
\end{abstract}
\maketitle



\section{Introduction}

The present paper deals with Poisson reduction of a cotangent bundle 
$T^* Q$ with respect to a Hamiltonian action by a compact Lie group 
$K$ that comes as the cotangent lifted action from the configuration 
manifold $Q$. 
We assume that $Q$ is Riemannian and $K$ acts on $Q$ by isometries. 
Further, we suppose that $Q$ is of single isotropy type, i.e., 
$Q = Q_{(H)}$ for some subgroup $H$ of $K$. 
Then the orbit space $Q/K$ is a smooth manifold.
However, in the presence of non-trivial isotropy, 
$H\neq\set{e}$, on the configuration 
space $Q$ one gets a non-trivial isotropy lattice on $T^*Q$,
whence $(T^*Q)/K$ cannot be a smooth manifold.
The cotangent bundle $T^*Q$ is equipped with its canonical symplectic 
form, and we have a standard momentum map $\mu : T^*Q \to \ko^*$. 

In Hochgerner \cite{H04} stratified symplectic reduction 
of $T^*Q$ was studied under these assumptions. In particular,
the following result 
was proved, following an approach that is generally called Weinstein
construction: 
Let $\orb$ be a coadjoint orbit lying in the image of the standard 
momentum map $\mu$. Then each smooth symplectic stratum $(T^*Q \sporb K)_{(L)}$ 
of the reduced space can be globally realized as  
\[
    (\WW\sporb K)_{(L)}   =   T^{*}(Q/K)\times_{Q/K}(\bscorb)_{(L)}/K 
\]  
where  
\[
  \WW :=  (Q\times_{Q/K}T^{*}(Q/K))\times_{Q}\bsc_{q\in Q}\ann{\ko_{q}}          
 \cong  T^{*}Q 
\] 
as symplectic manifolds with a Hamiltonian $K$-action. 
Moreover, the reduced symplectic structure in terms intrinsic to this 
realization was computed. These results were applied to Calogero-Moser 
systems with spin associated to polar representations of compact Lie 
groups.

In the above setting the reduced space $T^*Q/K$ is a (singular) 
Poisson space in a natural way. The underlying theorem due to Ortega 
and Ratiu \cite{OR98} is presented in \ref{sing-P-red}. The goal of 
this paper is to determine the reduced Poisson bracket on $T^*Q/K$.

In Section \ref{sing-geo} we give the necessary background from singular 
geometry.

Section \ref{W-bundle} presents the basic setting and the prerequisites 
of the paper, in particular, recalling from Hochgerner \cite{H04} the
construction of the Weinstein space $\WW$ and the intrinsic symplectic 
structure on it.

Finally, in Section \ref{g-P-red} we compute the reduced gauged Poisson 
bracket on the Weinstein realization $\WW/K$ of $T^*Q/K$. 
This is done by determining the Poisson bracket on $\WW$ via the 
symplectic structure intrinsic to this space and then using the 
identification $C^{\infty}(\WW/K) \cong C^{\infty}(\WW)^K$.
In particular, it will be important to introduce a 
suitably chosen linear 
connection $\widetilde{A}$ on the bundle 
\[
\WW \longto Q \times_{Q/K} T^*(Q/K).
\]
The formula for the reduced Poisson bracket on $\WW/K$ will involve 
the canonical Poisson structure on $T^*(Q/K)$, the pairing of a curvature 
term associated to the mechanical connection $A$ on $Q \to Q/K$ 
(see Section \ref{W-bundle}) with the appropriate 
point of $\bsc_{q \in Q} \ann\ko_q$, and the Lie bracket of the fiber 
derivatives on $\bsc_{q \in Q} \ann\ko_q$.

Note that this result generalizes the formula for the reduced Poisson
bracket on $T^*Q/K$ presented in Zaalani \cite{Zaa99} and in Perlmutter 
and Ratiu \cite{PR04} for the case that $K$ acts on $Q$ freely.

In Section \ref{symp-leaves} we compute the symplectic leaves of the 
reduced space $\WW/K$, and in Section \ref{sec:ex} we give 
examples of our constructions.

\section{Conventions}\label{s:transfgr} 
Let $K$ be a Lie group acting on a
manifold $M$.
In fact, we will only be concerned with the case where $K$ is compact,
$M$ is Riemannian, and $K$ acts on $M$ through isometries,
i.e. $M$ is a Riemannian $K$-space. 
The action will be written as $l: K\times M\to M$, $(k,x)\mapsto
l(k,x)=l_{k}(x)=l^{x}(k)=k.x$. 
Sometimes the action will be lifted to the tangent bundle $TM$. That
is, we will consider
$h.(x,v):=(h.x,h.v):=Tl_{h}.(x,v)=(l_{h}(x),T_{x}l_{h}.v)$ 
where $(x,v)\in TM$. 
As the
action is a transformation by a diffeomorphism it may also be lifted to
the cotangent bundle. This is the cotangent lifted action which is
defined by
$h.(x,p):=(h.x,h.p):=T^{*}l_{h}.(x,p)=
(h.x,T_{h.x}^{*}l_{h^{-1}}.p)$ where $(x,p)\in T^{*}M$. 

The \caps{fundamental vector field} is going to repeatedly
play an important role. It is defined by
\[
 \zeta_{X}(x) 
 := \dd{t}{}|_{0}l(\exp(+tX),x)
  = T_{e}l^{x}(X)
\]
where $X\in\ko$. The fundamental vector field mapping $\ko\to\X(M)$,
$X\mapsto\zeta_{X}$ is a Lie algebra anti-homomorphism. 
By definition the flow of
$\zeta_{X}$ is given by $l_{\exp(tX)}$. 

If the action by $K$ on $M$ is proper, in the sense that $K\times M\to
M\times M$, $(k,x)\mapsto(x,l(k,x))$ is a proper mapping, then we have
the Slice and Tube Theorem at our disposal. 
Properness is automatic for compact $K$.
An exposition of these
facts can be found in Palais and Terng \cite{PT88} or Ortega and Ratiu
\cite{OR04}, for example. 

Let $H$ be a subgroup of $K$.
A point $x \in M$ is said to be of \caps{isotropy} or \caps{orbit type} 
$H$ if its isotropy group $K_x = \{k \in K : k.x=x\}$ is conjugate to 
$H$ within $K$ for which we shall write $K_x \sim H$. 
The family of subgroups of $K$ conjugate to $H$ within $K$ is denoted by 
$(H)$ and called the \caps{conjugacy class} of $H$. 
We will deal with the \caps{isotropy} or \caps{orbit type submanifold} 
$M_{(H)} := \{x \in M : K_x \sim H\}$ of type $H$, the set 
$M_H := \{x \in M : K_x = H\}$ of points that have \caps{symmetry type} $H$, 
and the set $\fix{H} := M^H := \{x \in M : H \subseteq K_x\}$ of points that 
are fixed by $H$. Then $M_{(H)}$ is a submanifold of $M$, $\fix{H}$ is a
totally geodesic submanifold of $M$, and $M_H$ is an open submanifold of 
$\fix{H}$.

If $(M,\om)$ is a symplectic manifold we define the associated Poisson
bracket and Hamiltonian \vf by
\[
 \bra{f,g}
 =
 \om(\ham{f},\ham{g})
 = 
 -\ham{f}(g)
\]
where $f,g\in\cinf(M)$. This choice of sign is compatible with that in
\cite{OR04}. It is, however, not universal.

\section{Singular geometry} \label{sing-geo}

\subsection{Singular spaces and smooth structures}\label{sec:sg_smooth}
First we introduce the Whitney condition $(b)$ which will be
necessary in the definition of Whitney stratified spaces -- see
Definition \ref{def:stratif_space}. 
We follow the approach of Mather \cite{Mat70}.

\begin{definition}[Whitney condition $(b)$ in $\R^{n}$]
Let $X$, $Y$ be disjoint sub-manifolds of $\R^{n}$ with $\dim X=r$. The pair
$(X,Y)$ is said to satisfy \caps{condition} $(b)$ \caps{at} $y\in Y$
if the following is true. Consider sequences $(x_{i})_{i}$,
$(y_{i})_{i}$ in $X$, $Y$, respectively, such that $x_{i}\to y$ and
$y_{i}\to y$. Assume that $T_{x_{i}}X$
converges to some $r$-plane $\tau\subeq T_{y}\R^{n}=\R^{n}$, and
that the lines spanned by the vectors $y_{i}-x_{i}$ converge 
-- in $\R P^{n-1}$ --
to some line
$l\subeq\R^{n}=T_{y}\R^{n}$. Then $l\subeq\tau$. The pair $(X,Y)$
satisfies \caps{condition} $(b)$ if it does so at every $y\in Y$. 
\end{definition}

Obviously condition $(b)$ behaves well under diffeomorphisms in the
following sense: for $i=1,2$ consider pairs $(X_{i},Y_{i})$ in
$\R^{n}$, points $y_{i}\in Y_{i}$, open neighborhoods $U_{i}\subeq\R^{n}$ of
$y_{i}$, and a \diffeo $\phi: U_{1}\to U_{2}$ sending $y_{1}$ to
$y_{2}$ and satisfying $\phi(U_{1}\cap X_{1}) = U_{2}\cap X_{2}$ as
well as $\phi(U_{1}\cap Y_{1})=U_{2}\cap Y_{2}$. 
Thus it makes sense to formulate this condition for manifolds.

\begin{definition}[Whitney condition $(b)$]\label{def:Wit_b}
Let $M$
be a \mf and $X,Y$ disjoint sub-manifolds. Now $(X,Y)$ is said to satisfy
\caps{condition} $(b)$ if the following holds for all $y\in Y$. Let
$(U,\phi)$ be a chart around $y$. Then the pair 
$(\phi(X\cap U),\phi(Y\cap U))$ 
satisfies condition $(b)$ at $\phi(y)$. 
\end{definition}

By the above this definition is independent of the chosen chart in the
formulation.

\begin{ex}
Consider $M=\C^{3}=\set{(x,y,z)}$ with $Y$ the $z$-axis, and
$X=\set{(x,y,z): y^{2}+x^{3}-z^{2}x^{2}=0}\setminus Y$. Then the pair
$(X,Y)$ satisfies condition $(b)$ at
all points in $Y$ except at $y=0$. 
Notice that we can refine the decomposition of $M$ as
$(X,Y\setminus\set{0},\set{0})$. Now all pairs in this finer
decomposition satisfy condition $(b)$. 
\end{ex}

Let $X$ be a para-compact and second countable topological Hausdorff
space,  
and let $(I,\le)$ be a partially ordered set. 

\begin{definition}[Decomposed space]
An $I$-\caps{decomposition} of
$X$ is a locally finite partition of $X$ into smooth manifolds $S_{i}$, $i\in
I$ which are disjoint (but may consist of finitely many 
connected components with differing dimension), and satisfy:
\begin{enumerate}[\up (i)]
\item Each $S_{i}$ is locally closed in $X$;
\item $X=\bigcup_{i\in I}S_{i}$;
\item $S_{j}\cap\ov{S_{i}} \neq \emptyset \iff S_{j}\subeq\ov{S_{i}}
  \iff j\le i$. 
\end{enumerate}
The third condition is called \caps{condition of the frontier}. The
manifolds $S_{i}$ are called \caps{strata} or \caps{pieces}. In the
case that $j<i$ one often writes $S_{j}<S_{i}$ and calls $S_{j}$
\caps{incident} to $S_{i}$ or says $S_{j}$ is a \caps{boundary piece}
of $S_{i}$. 
\end{definition}

We define the dimension of a manifold consisting of finitely many
connected components to be the maximum of the dimensions of the
manifold's components. 

The \caps{dimension} of the decomposed space $X$ is defined as
\[
 \dim X :=\sup_{i\in I}\dim S_{i}
\]
and we will only be concerned with spaces where this supremum is
attained. 

The \caps{depth of the stratum} $S_{i}$ of the decomposed space $X$ 
is defined as
\begin{multline*}
 \dep S_{i} := \sup\{ l\in\N: \text{ there are strata }
   S_{i_{0}}=S_{i},S_{i_{1}},\dots,S_{i_{l}}\\
   \text{ such that } S_{i_{0}}<\ldots<S_{i_{l}} \},
\end{multline*}
Notice that $\dep S_{i}$ is always finite; indeed, else there would be
an infinite family $(S_{j})_{j\in J}$ with $S_{j}>S_{i}$ thus making
any neighborhood of any point in $S_{i}$ meet all of the $S_{j}$ which
contradicts local finiteness of the decomposition.
The \caps{depth} of $X$ is 
\[
 \dep X := \sup\set{\dep S_{i}: i\in I}.
\]
Thus, if $X$ consists of just one stratum, then $\dep X = 0$. From the
frontier condition we have that $\dep S_{i} \le \dim X - \dim S_{i}$,
and also $\dep X\le \dim X$. 

A simple example for a decomposed space is a \mf with boundary with
big stratum the interior and small stratum the boundary. 
Also manifolds with corners are decomposed spaces in the obvious way. 
Likewise the cone
$CM:=(M\times[0,\infty))/(M\times\set{0})$ over a manifold $M$ is a
decomposed space, the partition being that into cusp and open cylinder
$M\times(0,\infty)$. 

The following definition of singular charts and smooth structures on
singular spaces is due to Pflaum \cite[Section 2]{Pfl01}. 

\begin{definition}[Singular charts]
Let $X=\bigcup_{i\in I}S_{i}$ be a decomposed space. A
\caps{singular chart} $(U,\psi)$ with patch $U$ an open subset of $X$
is to satisfy the following. 
\begin{enumerate}[\up (i)]
\item 
$\psi(U)$ is locally closed in $\R^{n}$; 
\item 
$\psi: U\to\psi(U)$ is a homeomorphism;
\item 
For every stratum $S_{i}$ that meets $U$ the restriction
$\psi|S_{i}\cap U: S_{i}\cap U\to\psi(S_{i}\cap U)$ is a
diffeomorphism onto a smooth sub-manifold of $\R^{n}$. 
\end{enumerate}
Two singular charts $\psi: U\to\R^{n}$ and $\phi: V\to\R^{m}$ are
called \caps{compatible at} $x\in U\cap V$ if there is an open
neighborhood $W$ of $x$ in $U\cap V$, a number $N\ge\max\set{n,m}$, and a
\diffeo $f: W_{1}\to W_{2}$ between open subsets of $\R^{N}$ such that: 
\[
\xymatrix{
 & W \ar[d]_{\psi}
     \ar[dr]^{\phi}\\
 & \psi(W) 
     \ar[r]_{f|\psi(W)}
     \ar @{^{(}->}[d]
  &\phi(W)
     \ar @{^{(}->}[d]\\
 {\R^{N}}
  &W_{1}
     \ar[r]^{f}
     \ar @{^{(}->}[l]
  &W_{2}
     \ar @{^{(}->}[r]
  &{\R^{N}}   
}
\]
It follows that $f|\psi(W):\psi(W)\to\phi(W)$ is a
homeomorphism. Further, for all strata $S$ that meet $W$ the
restriction $f|\psi(W\cap S): \psi(W\cap S)\to\phi(W\cap S)$ is a
\diffeo of sub-manifolds of $\R^{N}$. The charts $(U,\psi)$ and
$(V,\phi)$ are called \caps{compatible} if they are so at every point
of the intersection $U\cap V$. It is straightforward to check that
compatibility of charts defines an equivalence relation. 

A family of
compatible singular charts on $X$ such that the union of patches
covers all of $X$ is called a \caps{singular atlas}. Two singular
atlases are said to be compatible if all charts of the first are
compatible with all charts of the second. Again it is clear that
compatibility of  atlases forms an equivalence relation. 

Let $\mathfrak{A}$ be a singular atlas on $X$. Then we can consider the
family of all singular charts that belong to some atlas compatible
with $\mathfrak{A}$ to obtain a maximal atlas $\mathfrak{A}_{\textup{max}}$. 
\end{definition}

\begin{definition}[Smooth structure]\label{def:sg_smooth}
Let $X=\bigcup_{i\in I}S_{i}$ be a decomposed space. A maximal atlas
$\mathfrak{A}$ on $X$ is called a \caps{smooth structure} on the singular
space $X$. A continuous function $f:X\to\R$ is said to be
\caps{smooth} if the following holds. For all charts $\psi:
U\to\R^{n}$ of the atlas $\mathfrak{A}$ there is a smooth function $F:
\R^{n}\to\R$ such that $f|U=F\circ\psi$. The set of all smooth
functions on $X$ is denoted by $\cinf(X)$. 

A continuous map $f: X\to Y$ between decomposed spaces with smooth
structures is called \caps{smooth} if $f^{*}\cinf(Y)\subeq\cinf(X)$. 
An \caps{isomorphism} $F:X\to Y$ between decomposed spaces is a
homeomorphism that is smooth in both directions and maps strata of $X$
diffeomorphically onto strata of $Y$. 
\end{definition}

The smooth structure thus defined on decomposed spaces is in no way
intrinsic but is a structure that is additionally defined to do
analysis on decomposed spaces. 
Also a smooth map $f: X\to Y$ between decomposed spaces need not at all be
strata preserving. 

\begin{definition}[Cone space]
A decomposed space $X=\bigcup_{i\in I}S_{i}$ is called a \caps{cone
space} if the following is true. Let $x_{0}\in X$ arbitrary and $S$ the
stratum passing through $x_{0}$. Then there is an open neighborhood $U$ of
$x_{0}$ in $X$, there is a decomposed space $L$ with global
chart $\psi: L\to S^{l-1}\subeq\R^{l}$, and furthermore there is an
\iso of decomposed spaces
\[ 
 F: U\to(U\cap S)\times CL
\] 
such that $F(x)=(x,c)$ for all $x\in U\cap S$. Here
$CL=(L\times[0,\infty))/(L\times\set{0})$ is decomposed into the cusp
$c$ on the one hand, while the other pieces are of the form stratum of
$L$ times $(0,\infty)$. Thus we can take 
 $\Psi: CL\to\R^{l}$,
 $[(z,t)]\mapsto t\psi(z)$ 
as a global chart on $CL$ thereby defining
a smooth structure on $CL$ whence also on the product 
$(U\cap S)\times CL$. 

The space $L$ is called a \caps{link}, and the chart $F$ is referred
to as a \caps{cone chart} or also \caps{link chart}. Of course, the
link $L$ depends on the chosen point $x_{0}\in X$.
\end{definition}
 
%

An example for a cone space is the quadrant
$
 Q := \set{(x,y)\in\R^{2}: x\ge0\text{ and }y\ge0}
$.
A typical
neighborhood of $0\in Q$ is of the form 
$
 \set{(x,y): 0\le x<r\text{ and }0\le y<r}
$.
The link with respect to the point $0$ then is
the arc 
$
 L := \set{(\cos\var,\sin\var): 0\le\var\le\by{\pi}{2}}
$.
More generally manifolds with corners carry the structure of cone
spaces.

\begin{definition}[Stratified spaces]\label{def:stratif_space}
Let $X\subeq\R^{m}$ be a subset and assume that $X$ is a decomposed
space, i.e. $X=\bigcup_{i\in I}S_{i}$, and that the strata $S_{i}$ be
sub-manifolds of $\R^{m}$. 
The $I$-decomposed space $X$ is said to be \caps{(Whitney) stratified}
if all pairs $(S_{i},S_{j})$ with $i>j$ satisfy condition $(b)$ -- see
Definition \ref{def:Wit_b}.
For sake of convenience we will simply say stratified instead of Whitney
stratified. 
\end{definition}

\begin{theorem}\label{thm:stratif_vs_cone}
Let $X\subeq\R^{m}$ be a subset of a Euclidean space and assume that
$X=\bigcup_{i\in I}S_{i}$ is decomposed. Then $X$ is stratified if and
only if $X$ is a cone space.
\end{theorem} 

\begin{proof}
It is proved in Pflaum \cite{Pfl01a} that every (Whitney) stratified
space is also a cone space. 

An outline of the converse direction is given in Sjamaar and Lerman
\cite[Section 6]{SL91}, and also in Goresky and MacPherson \cite[Section 1.4]{GorMac88}. This
argument makes use of Mather's control theory as introduced in Mather
\cite{Mat70} as well as Thom's First Isotopy Lemma.
\end{proof}

The above theorem depends crucially on the fact that the decomposed
space $X$ can be
regarded as a subspace of some Euclidean space. As this assumption
will always be satisfied in the present context we will take the words
cone space and stratified space to be synonymous. In fact, Sjamaar
and Lerman \cite{SL91} take cone space to be the definition of
stratified space.

\begin{ex}\label{ex:quot_smooth}
As an example consider a compact Lie group $K$
acting by isometries on a smooth Riemannian manifold
$M$. We are concerned with the orbit projection 
$\pi: M\to M/K$ and
endow the orbit space with the final topology with respect to the 
projection map. 
For  basics on compact transformation groups see Bredon \cite{Bre72},
Palais and Terng \cite{PT88}, 
or Hochgerner \cite[Section 7]{H04}. 
Fix a point $x_{0}\in M$ with isotropy group
$K_{x_{0}}=H$. The slice representation is then the action by $H$ 
on $\nor_{x_{0}}(K.x_{0}) = T_{x_0}(K.x_0)^{\bot}$. 
By the Tube Theorem  there is a $K$-invariant open neighborhood $U$ 
of the orbit $K.x_{0}$ such that $K\times_{H}V\cong U$ 
as smooth $K$-spaces where $V$ is 
an $H$-invariant open neighborhood of $0$ in $\nor_{x_{0}}(K.x_{0})$.
 
Now let $p=(p_{1},\dots,p_{k})$ be a Hilbert basis for the algebra
$\poly(V)^{H}$ of $H$-invariant polynomials on $V$. That is,
$p_{1},\dots,p_{k}$ is a finite system of generators
for $\poly(V)^{H}$. The Theorem of Schwarz \cite[Theorem 1]{Schwarz75}
now says that $p^{*}:
\cinf(\R^{k})\to\cinf(V)^{H}$ is surjective. 
Moreover, the induced mapping $q: V/H\to\R^{k}$
is continuous, injective, and proper. See also Michor \cite{Mic97}. 

Consider the isotropy type
sub-manifolds $M_{(H)}$
which is the manifold of all points of $M$ whose isotropy subgroup is
conjugate to $H$ within $K$.  
These give a $K$-invariant decomposition of
$M$ as $M = \bigcup_{(H)}M_{(H)}$ where $(H)$ runs through the
isotropy lattice of the $K$-action on $M$. We thus get a decomposition
of the orbit space 
\[
 M/K 
 =
 \bigcup_{(H)}M_{(H)}/K
\]
where again $(H)$ runs through the isotropy lattice of the $K$-action
on $M$. It is well-known (e.g. Palais and Terng \cite{PT88}) 
that this decomposition renders $M/K$ a decomposed space. 

Now a theorem of Pflaum \cite[Theorem 5.9]{Pfl01} says that the
induced mapping $\psi: U/K\to\R^{k}$ as defined in  the diagram 
\[
\xymatrix{
 U\ar[r]^-{{\simeq}}\ar @{->>}[d]
 &{K\times_{H}V}\ar @{->>}[d]\\
 {U/K}\ar[r]^{{\simeq}}_{{\phi}}\ar @{-->}[dr]_{{\psi}}
 &{V/H}\ar[d]_{q}
 &V\ar @{->>}[l]\ar[dl]^{p}\\
 &{\R^{k}}
}
\]
is a typical singular chart around the point $K.x_{0}$ in the orbit
space. Furthermore, the smooth functions with respect to this smooth
structure are $\cinf(M/K)=\cinf(M)^{K}$, i.e. none other than the
$K$-invariant smooth functions on $M$: indeed, by Schwarz' Theorem we
have
\[
 \psi^{*}\cinf(\R^{k})
 = 
 \phi^{*}q^{*}\cinf(\R^{k})
 =
 \phi^{*}\cinf(V)^{H}
 =
 \cinf(U)^{K}
\]
whence $\cinf(U)^{K}=\cinf(U/K)$. Finally, the decomposition of $M/K$
by orbit types turns the orbit space into a stratified space with
smooth structure. 

Note that Bierstone \cite{Bier75} showed that the semi-analytic 
stratification of the orbit space $p(V)$ of a linear $K$-space $V$ 
coincides with its stratification by components of sub-manifolds of 
given isotropy type. Here $K$ is a compact Lie group and $p$ is a 
Hilbert basis as above. Further, the semi-analytic stratification 
of $p(V)$ satisfies Whitney's condition $(b)$.
\end{ex}

\subsection{Singular Poisson reduction}\label{sing-P-red}

Let $K$ be a Lie group acting properly on a smooth manifold $M$.
We equip the orbit space $M/K$ with the quotient topology with respect
to the canonical projection $\pi : M \to M/K$.
The set of smooth functions on $M/K$ is defined by the requirement that 
$\pi$ is a smooth map, i.e.,
\[
C^{\infty}(M/K) := \{f \in C^0(M/K) : f \circ \pi \in C^{\infty}(M)\}.
\]

\begin{Theorem}[Singular Poisson reduction]
Let $(M,\{ \cdot,\cdot \})$ be a Poisson manifold, $K$ a Lie group, and let 
$l : K \times M \to M$ be a smooth proper Poisson action, i.e., 
$l_k^*\{f,g\} = \{l_k^*f,l_k^*g\}$
for 
$f,g \in C^{\infty}(M)$ and $k \in K$. 
Then we have: 
\begin{itemize}
\item[(i)] The pair $(C^{\infty}(M/K),\{ \cdot,\cdot \}^{M/K})$ is a
           Poisson algebra, where the Poisson bracket 
	   $\{ \cdot,\cdot \}^{M/K}$ is characterized by 
	   $\{f,g\}^{M/K} \circ \pi = \{f \circ \pi,g \circ \pi\}$, 
	   for any $f,g \in C^{\infty}(M/K)$, and $\pi : M \to M/K$
	   denotes the canonical smooth projection.
\item[(ii)] Let $h \in C^{\infty}(M)^K$ be a $K$-invariant function on $M$.
           The flow $\fl_t$ of the Hamiltonian vector field
	   associated to $h$
	   commutes with the $K$-action, so it induces a flow $\fl_t^{M/K}$
	   on $M/K$ which is Poisson and is characterized by the identity
	   $\pi \circ \fl_t = \fl_t^{M/K} \circ \pi$.
\item[(iii)] The flow $\fl_t^{M/K}$ is the unique Hamiltonian flow defined 
           by the function $H \in C^{\infty}(M/K)$ which is given by 
           $H \circ \pi =h$. 
\end{itemize}
\end{Theorem}

\proof
This is due to Ortega and Ratiu \cite{OR98}.
\endproof

If, in particular, $K$ is a compact Lie group acting by isometries on a 
smooth Riemannian manifold $M$, then by the Tube Theorem \cite{PT88} and 
Schwarz' Theorem \cite{Schwarz75} we may identify $C^{\infty}(M/K)$ with 
$C^{\infty}(M)^K$, 
the set of $K$-invariant smooth functions on $M$.
See Example \ref{ex:quot_smooth} for details.

\begin{definition}[Poisson stratified space]\label{def:poiss_stratif}
Let $X$ be a stratified space endowed with a smooth structure in the
sense of Subsection \ref{sec:sg_smooth}. Then $X$ is said to be a
\caps{singular Poisson space} if there is a Poisson bracket
\[
 \bra{\cdot,\cdot}:
 \cinf(X)\times\cinf(X)\longto\cinf(X)
\]
on the algebra of smooth functions determined by the smooth structure
such that the inclusion of each stratum $S\hookto X$ is a Poisson
morphism. In particular, the strata $S$ are Poisson manifolds in the
usual sense.
\end{definition}

An alternative definition of a singular Poisson space in terms of a
stratified Poisson bivector is given by 
Pflaum \cite{Pfl01}.

\begin{proposition}[Reduced Poisson structure]
Let $(M,\bra{\cdot,\cdot})$ be a Poisson manifold, $K$ a compact Lie group, and let 
$K$ act on $M$ by Poisson morphisms. 
Then $(M/K,\bra{\cdot,\cdot}^{M/K})$ is a singular Poisson
space.
\end{proposition}

\begin{proof}
By Example \ref{ex:quot_smooth}
the algebra $\cinf(M/K)$ is indeed determined by a smooth structure on
$M/K$. Thus it only remains to check that the inclusion of each
stratum $M_{(H)}/K\hookto M/K$ is a Poisson morphism. This is,
however, obvious.
\end{proof}

\subsection{Singular fiber bundles}\label{sec:sg_fb}

\begin{definition}[Singular fiber bundles]\label{def:sg_fb} 
Let $F$ and $P$ be stratified spaces (Definition \ref{def:stratif_space})
with smooth structure (Definition \ref{def:sg_smooth})
and $M$ be a smooth manifold. We
say that the topological fiber bundle 
\[
\xymatrix{
 {F}\ar @{^{(}->}[r]
 &{P}\ar[r]^-{\pi}
 &M
}
\]
is a \caps{singular fiber bundle}
if for each trivializing patch $U\subeq M$ the homeomorphism 
\[
 P|U\cong U\times F
\]
is an \iso of stratified spaces.
\end{definition}

There are two reasons for defining singular fiber bundles in this
way. Firstly, it is the kind of structure encountered in 
Hochgerner \cite[Theorem 5.5]{H04}, and secondly by Mather's control theory
\cite{Mat70} these bundles possess many features similar to ordinary
smooth fiber bundles. Some of these aspects are presented in the
subsections below.

Note that if $M$ is a Riemannian \mf which is acted upon by a compact
Lie group $K$ through isometries then the orbit projection mapping $M\toto
M/K$ is, in general, not a singular fiber bundle according to this
definition. Indeed, the fiber type of $M\toto M/K$ need not be locally
constant. 

\begin{lemma}\label{lem:sg_fb}
Let $\pi: P\to M$ be a singular fiber bundle with typical fiber
$F$. 
Let $S$ be a stratum of $P$. Then $\pi|S: S\to M$ is a smooth fiber
bundle.
\end{lemma}

\begin{proof}
Indeed, locally the stratum $S$ is diffeomorphic to $U\times S_{F}$
where $S_{F}$ is a stratum of $F$ and $U$ is a trivializing
neighborhood in $M$.
\end{proof}

\begin{definition}[Singular symplectic fiber bundles]\label{def:sg-sp-fb}
Let $F$ and $P$ be stratified symplectic spaces 
(i.e.\ singular Poisson spaces with smooth symplectic strata such
that the inclusion mappings are Poisson morphisms)
with smooth structure 
and $M$ be a smooth symplectic manifold. We
say that the singular fiber bundle 
\[
\xymatrix{
 {F}\ar @{^{(}->}[r]
 &{P}\ar[r]^-{\pi}
 &M
}
\]
is a \caps{singular symplectic fiber bundle}
if for each trivializing patch $U\subeq M$ the homeomorphism 
\[
 P|U\cong U\times F
\]
is an \iso of stratified symplectic spaces with respect to the
inherited symplectic structures. It follows, in particular, that $\pi$
is a Poisson morphism.
\end{definition}

\subsection{Control data}
The theory of control data is due to Mather \cite{Mat70}, and we
follow in our presentation of the subject that of \cite{Mat70}. 
Let $N$ be a smooth manifold, and $X\subeq N$ a stratified subset 
endowed with the relative topology
with
strata $S_{i}$ where $i\in I$ as in Subsection \ref{sec:sg_smooth}.
 
A \caps{tubular neighborhood} of a stratum $S_{i}$ in $X$ is a closed
neighborhood of $S_{i}$ in $N$ which is diffeomorphic to an inner product
bundle $\pi_{i}: E_{i}\to S_{i}$. Via the inner product we can measure the
vertical distance of a point in $E_{i}$ to $S_{i}$ and call this the
tubular neighborhood function $\rho_{i}:E_{i}\to\R$. Clearly, $\rho_{i}(x)=0$
if and only if $x\in S_{i}$. We can also think of the tubular neighborhood as being
retracted onto $S_{i}$ via the projection $\pi_{i}$. \caps{Control data}
associated to the stratification $\set{S_{i}: i\in I}$ 
of $X$ is a system of
tubular neighborhoods $\pi_{i}: E_{i}\to S_{i}$ satisfying the
following commutation relations:
\[
 (\pi_{j}\circ\pi_{i})(x) = \pi_{j}(x), \text{ and }
 (\rho_{j}\circ\pi_{i})(x) = \rho_{j}(x)
\]
whenever $j < i$ and both sides are defined.

\begin{proposition}\label{prop:controldata}
Suppose there exist control data to the stratification 
$\set{S_{i}: i\in I}$ of $X$. 
If $M$ is another manifold, and $f: N\to M$ a 
smooth mapping such that $f|S_{i}: S_{i}\to M$ is a submersion for all
$i\in I$ 
then the
control data may be chosen so that $f|S_i \circ \pi_{S_{i}}=f$ for all 
$i\in I$. 
\end{proposition}

\begin{proof}
See Mather \cite[Proposition 7.1]{Mat70}.
\end{proof}

If $f: N\to M$ is as in Proposition \ref{prop:controldata} then $f$ is
said to be a \caps{controlled submersion} from $X$ to $M$. 

By a \caps{stratified vector field} $\eta$ on $X$ we mean a collection
$\set{\eta_{i}:i\in I}$ where each $\eta_{i}$ is a smooth \vf on
$S_{i}$. Assume we are given a system of control data associated to
the stratification of $X$, and identify the tubular neighborhoods of
the strata with the corresponding inner product bundles. Then the
stratified \vf $\eta$ on $X$ is said to be a \caps{controlled vector
  field} if the following conditions are met.
For any stratum $S_{j}$ there is an open neighborhood $B_{j}$ of
$S_{j}$ in the tubular neighborhood $E_{j}$ such that for any 
stratum $S_{i}$ with $i>j$ the conditions
\[
 \ld_{\eta_{i}}(\rho_{j}|B_{j}\cap S_{i}) = 0,\text{ and }
 T_{x}(\pi_{j}|B_{j}\cap S_{i}).\eta_{i}(x)
   = \eta_{j}(\pi_{j}(x))
\]
are satisfied for all $x\in B_{j}\cap S_{i}$.

Let $J$ be an open neighborhood of $\set{0}\times X$ in $\R\times X$, and
assume that $\alpha: J\to X$ is a local one-parameter group 
which is smooth in the sense of Definition
\ref{def:sg_smooth}. We say that $\alpha$ \caps{generates} the
stratified \vf $\eta$ if $J$ is maximal such that 
each stratum $S_{i}$ is invariant under $\alpha$ and 
$
 \dd{t}{}|_{0}\alpha(t,x) = \eta_{i}(x)
$
for all
$x\in S_{i}$ and all $i\in I$.  

\begin{proposition}
Assume $\eta$ is a controlled \vf on $X$. Then there is a unique
smooth one-parameter group which generates $\eta$.
\end{proposition}

\begin{proof}
See Mather \cite[Proposition 10.1]{Mat70}.
\end{proof}

\begin{proposition}
Assume $f:N\to M$ is a smooth map such that $f|X: X\to M$ is a
controlled surjective submersion.
Then the following are true.
\begin{itemize}
\item 
Let $\xi$ be a smooth \vf on $M$. Then
there is a controlled \vf $\eta$ on $X$ such that $\eta_{i}$ and $\xi$
are $f|S_{i}$-related for all $i\in I$. 
\item
Suppose further that $f|X : X \to M$ is a proper map. Then $f|X : X\to M$ is a
singular fiber bundle.
\end{itemize}
\end{proposition}

\begin{proof}
See Mather \cite[Proposition 9.1]{Mat70} for the first statement. 
Concerning the second assertion, \cite[Proposition 11.1]{Mat70} 
states that under these assumptions the mapping
$f|X : X\to M$ is locally topologically trivial, and it follows from
\cite[Corollary 10.3]{Mat70} that the trivializing homeomorphisms are,
in fact, isomorphisms of stratified spaces. Thus $f|X : X\to M$ is a
singular fiber bundle in the sense of Definition \ref{def:sg_fb}.
\end{proof}

\subsection{Pullback bundles} \label{sec:pullbacks}
Let $M$ and $Y$ be smooth manifolds, and let $\tau: Y\to M$ be a
smooth mapping. Consider a singular fiber bundle $\pi: X\to M$ with
typical fiber $F$ as in Definition \ref{def:sg_fb}. We consider
further the topological pullback bundle of $X$ and $Y$ over $M$ with
the following notation.
\[
\xymatrix{
 {X\times_{M}Y}\ar @{=}[r]&
 {\tau^{*}X}\ar[r]^-{\pi^{*}\tau}\ar[d]_-{\tau^{*}\pi}&
 {X}\ar[d]^-{\pi}\\ 
 &{Y}\ar[r]^-{\tau}&{M}
}
\]
Now we can endow $X\times_{M}Y$ with the product stratification given
by strata of the form $S\times_{M}Y$ which is the smooth fibered
product of a stratum $S$ of $X$ with $Y$ over $M$. Note that
$S\times_{M}Y$ is a well defined pull back bundle by Lemma \ref{lem:sg_fb}.
Moreover, $X\times_{M}Y$ inherits a smooth structure in the sense of
Subsection \ref{sec:sg_smooth} from the canonical topological inclusion
$X\times_{M}Y\hookto X\times Y$. The singular space with smooth
structure thus obtained is called the \caps{fibered product} of $X$
and $Y$ over $M$. Since $\pi: X\to M$ is a singular fiber bundle it
follows that 
$
 \tau^{*}\pi: X\times_{M}Y\to Y
$
is a singular fiber
bundle as well with the same typical fiber $F$.
Moreover, this construction satisfies the following universal
property. Let $Z$ be a singular space with smooth structure and
$f_{1}: Z\to X$, $f_{2}: Z\to Y$ be smooth mappings satisfying
$\pi\circ f_{1}=\tau\circ f_{2}$. Then there is a unique smooth map
$f=(f_{1},f_{2})$ such that the following commutes.
\[
\xymatrix{
 &{Z}\ar[dl]_-{f_{1}}\ar @{-->}[d]_-{f}
     \ar[dr]^-{f_{2}}\\
 {X}\ar[dr]^-{\pi}&
 {X\times_{M}Y}\ar[l]_-{\pi^{*}\tau}\ar[r]^-{\tau^{*}\pi}&
 {Y}\ar[dl]_-{\tau}\\
 &{M}
}
\]
Therefore, in this sense, pull backs exist in the category of singular
fiber bundles. There is, in fact, a similar notion of pull backs in
the work of Davis \cite{Dav78}.

\section{Mechanical connection and Weinstein construction} \label{W-bundle}

Suppose that $Q$ is a Riemannian manifold, and $K$ is a compact Lie
group which acts on $Q$ by isometries. Moreover, $Q$ is assumed to be
of single isotropy type, i.e. $Q=Q_{(H)}$ where $H$ is an isotropy
subgroup of $K$. The $K$-action then induces a Hamiltonian action
on the cotangent bundle $T^{*}Q$ by cotangent lifts. This means that
the lifted action respects the canonical symplectic form
$\Om=-d\theta$ on $T^{*}Q$ and there is a 
\momap $\mu: T^{*}Q\to\ko^{*}$ 
given by 
$\vv<\mu(q,p),X>
 = \theta(\zeta^{T^{*}Q}_{X})(q,p)
 = \vv<p,\zeta_{X}(q)>$ 
where
$(q,p)\in T^{*}Q$, $X\in\ko$, $\zeta_{X}$ is the fundamental
vector field associated to the $K$-action on $Q$, 
and $\zeta^{T^{*}Q}_{X}\in\X(T^{*}Q)$ is the
fundamental vector field associated to the cotangent lifted action. 

\subsection{Mechanical connection}
Since the $K$-action on $Q$ has only a single isotropy type, the orbit
space $Q/K$ is a smooth manifold, and the projection $\pi: Q\to Q/K$
is a surjective Riemannian submersion with compact fibers.
However, the lifted action by $K$ on
$T^{*}Q$ is already much more complicated, and the quotient space $(T^{*}Q)/K$
is only a stratified space in general. Its strata are of the form
$(T^{*}Q)_{(L)}/K$ where $(L)$ is in the isotropy lattice of $T^{*}Q$. 

The vertical sub-bundle of $TQ$ with respect to $\pi: Q\to Q/K$ is
$\ver := \ker T\pi$. Via the $K$-invariant Riemannian metric we obtain
the horizontal sub-bundle as $\hor := \ver^{\bot}$. We define the dual
horizontal sub-bundle of $T^{*}Q$ as the sub-bundle $\hor^{*}$
consisting of those co-vectors that vanish on all vertical
vectors. Likewise, we define the dual
vertical sub-bundle of $T^{*}Q$ as the sub-bundle $\ver^{*}$
consisting of those co-vectors that vanish on all horizontal
vectors.

We choose and fix a $K$-invariant inner product on $\ko$. 
For $X,Y\in\ko$ and $q\in Q$ we define
$\ine_{q}(X,Y) := \vv<\zeta_{X}(q),\zeta_{Y}(q)>$ and call this  
the \caps{inertia tensor}.  
This gives a
non-degenerate pairing on $\ko_{q}^{\bot}\times\ko_{q}^{\bot}$, whence
it gives an identification $\check{\ine_{q}}:
\ko_{q}^{\bot}\to(\ko_{q}^{\bot})^{*} = \ann\ko_{q}$. We use this \iso
to define a one-form on $Q$ with values in the bundle 
$\bsc_{q\in Q}\ko_{q}^{\bot}$ by the following: 
\[
\xymatrix{
 {T_{q}^{*}Q} \ar[r]^-{\mu_{q}} 
              &
 {\ann\ko_{q}}
 \ar[d]^{(\check{\mathbb{I}_{q}})^{-1}}
             \\
 {T_{q}Q}     \ar[u]^{\simeq}
              \ar @{-->}[r]^{A_{q}}& 
 {\ko_{q}^{\bot}}
}
\]
Here the isomorphism $T_q Q \stackrel{\simeq}{\to} T^*_q Q$ is obtained 
via the $K$-invariant Riemannian metric on $Q$. 
See Smale \cite{Sma70} or
Marsden, Montgomery, and Ratiu \cite[Section 2]{MMR90}. 
The form $A$ shall be called the \caps{mechanical connection} on $Q\to
Q/K$. It has the following properties.
It follows from its definition that 
$TQ\to\bsc_{q\in Q}\ko_{q}^{\bot}$, 
$(q,v)\mapsto(q,A_{q}(v))$
is equivariant, 
$\ker A_{q}=T_{q}(K.q)^{\bot}$, 
and 
$A_{q}(\zeta_{X}(q))=X$ for all $X\in\ko_{q}^{\bot}$. 

This means that $A\in\Om^{1}(Q;\ko)$ given by 
$A: TQ\to\ko_{q}^{\bot}\hookto\ko$, 
$(q,v)\mapsto A_{q}(v)$
is a principal connection form on the $K$-manifold $Q$ in the sense of
Alekseevsky and Michor \cite[Section 3.1]{AM95}. According to
\cite[Section 4.6]{AM95} the curvature form associated to $A$ is
defined by 
\[
 \curv^{A}
 :=
 dA-\by{1}{2}[A,A]^{\wedge}
\]
where 
\[
 [\var,\psi]^{\wedge} (v_1,\ldots,v_{l+k})
 :=
 \by{1}{k!l!}
  \sum_{\sigma}\sign\sigma
  [\var(v_{\sigma1},\dots,v_{\sigma l}),\psi(v_{\sigma(l+1)},\dots,v_{\sigma(l+k)})]
\]
is the graded Lie bracket on 
$
 \Om(Q;\ko) 
 := \bigoplus_{k=0}^{\infty}\Gamma(\Lam^{k}T^{*}Q\otimes\ko)
$, 
and $\var\in\Om^{l}(Q;\ko)$
and $\psi\in\Om^{k}(Q;\ko)$. 
The minus in the definition of $\curv^A$ emerges, since we are dealing with 
left actions. 
The form $\curv^A$ is called the \caps{mechanical curvature}.

\begin{proposition}[Properties of $\curv^A$]\label{prop:curvA}
Let $Z_1,Z_2\in\ko$ and $v,w$ horizontal vector fields on $Q$ with
respect to the orbit projection $\pi: Q\toto Q/K$. Then:
\begin{enumerate}[\up (i)]
\item
$\curv^A(\zeta_{Z_{1}},\zeta_{Z_{2}})(q)\in\ko_q$.
\item
$\curv^A(v,\zeta_{Z_{2}})(q) = 0$.
\item
$\curv^A(v,w)(q)\in\ko_{q}^{\bot}$.
\end{enumerate}
Moreover, the form $\curv^A: \Lam^{2}TQ\to\ko$ is $K$-equivariant,
$\Ad(h).\curv^A(v,w)(q) = \curv^A(v,w)(q)$ for all $h\in K_q$, 
and $\curv^A$ drops to a
well defined form 
\[
 \curv^{A}_{0}\in\Om^2(Q/K;(\bsc_{q\in Q}\fix{K_q}\cap\ko_q^{\bot})/K)
\] 
where $\fix{K_q} := \set{X\in\ko: K_q \subeq K_X}$. 
\end{proposition}

Notice that 
$(\bsc_{q\in Q}\fix{K_q}\cap\ko_q^{\bot})/K$ is a smooth manifold and
a topological subspace of the stratified space 
$(\bsc_{q\in Q}\ko_q^{\bot})/K$ -- see Lemma \ref{E:Q}.

\begin{proof}
Assertion (i) is true since the vertical bundle is integrable: Indeed,
$\Phi := \zeta\circ A: TQ\to\ver$ defines a (generalized) principal bundle
connection on $Q\toto Q/K$ in the sense of Alekseevsky and Michor
\cite{AM95}. As usual, the curvature $R$ associated to $\Phi$ is given
by $R(X,Y) = \Phi[X-\Phi X,Y-\Phi Y]$. By \cite[Proposition 4.4]{AM95}
$R$ and $\curv^A$ are related by 
\[
 R = -\zeta\circ\curv^A
\]
whence (i) follows. Using this relation again it follows that 
$\curv^A(v,\zeta_{Z_{2}})(q)\in\ko_{q}$ and this is already sufficient
for the purpose of this paper. However, using \cite[Proposition
  4.7]{AM95} and the Slice theorem for Riemannian actions the stronger
result (ii) is true as well. Assertion (iii) follows also by using the
Slice theorem. 

Clearly $\curv^A$ is $K$-equivariant. The element $\curv^A(v,w)(q)$ is
fixed by all $h\in K_q$ since this is true for $q$ and horizontal
vectors at $q$ have isotropy at least $K_q$. Finally, the last
conclusion is an obvious consequence of the above. It is only to notice
that $(\bsc_{q\in Q}\fix{K_q}\cap\ko_q^{\bot})/K$ really is a smooth
fiber bundle over $Q/K$ which, again, is a consequence of the Slice
theorem. 
\end{proof}

We define a point-wise dual 
$A_{q}^{*}: 
 \ann\ko_{q}\to
 \ver_{q}^{*}\subeq T_{q}^{*}Q$ 
by the formula 
$
 A_{q}^{*}(\lam)(v) 
 = 
 \lam(A_{q}(v))
$
where $\lam\in\ann\ko_{q}$ and $v\in T_{q}Q$. Notice that 
$
 A_{q}^{*}(\mu_{q}(p)) = p
$
for all $p\in\ver_{q}^{*}$
and 
$
 \mu_{q}(A_{q}^{*}(\lam)) = \lam
$
for all $\lam\in\ann\ko_{q}$.

\subsection{Weinstein bundle construction}
Using the horizontal lift
mapping which identifies $\hor \cong
(Q\times_{Q/K}T(Q/K))$ and the mechanical connection
$A$ 
we obtain an \iso 
\[
 TQ = \hor\oplus\ver 
 \longto
  (Q\times_{Q/K}T(Q/K))\times_{Q}\bsc_{q\in Q}\ko_{q}^{\bot}
\]
of bundles over $Q$. Via the Riemannian structure there is a dual 
version to this \iso, and to save on typing we will abbreviate 
\[
 \WW :=
 (Q\times_{Q/K}T^{*}(Q/K))\times_{Q}\bsc_{q\in Q}\ann{\ko_{q}}
 \cong
 \hor^{*}\oplus\ver^{*}.
\]
To set up some notation for the upcoming proposition, 
and clarify the picture consider the following
stacking of pull-back diagrams.
It will be reference point for the whole paper.
\begin{equation} \label{W-diag}
\xymatrix{
 {\WW} \ar[r]^-{\rho^{*}\widetilde{\tau}=\widetilde{\widetilde{\tau}}}
       \ar[d]_{\widetilde{\tau}^{*}\rho=\widetilde{\rho}} &
 {\bsc_{q}\ann\ko_{q}}
       \ar[d]^{\rho} \\ 
 {Q\times_{Q/K}T^{*}(Q/K)}
       \ar[r]^-{\pi^{*}\tau=\widetilde{\tau}}
       \ar[d]_{\tau^{*}\pi=\widetilde{\pi}} &
 Q 
       \ar[d]^{\pi}  \\
 T^{*}(Q/K)  
       \ar[r]^-{\tau} &
 Q/K       
}
\end{equation}
The upper stars in this diagram are, of course, not pull-back
stars. It is in fact the transition functions that are being pulled-back,
whence the name.     

\begin{proposition}[Symplectic structure on \WW]\label{prop:WW}
There is a dual \iso
\begin{align*}
 \psi=\psi(A):\; &
 (Q\times_{Q/K}T^{*}(Q/K))\times_{Q}\bsc_{q\in Q}\ann\ko_{q} = \WW
  \longto T^{*}Q, \\
 &(q,\eta,\lambda)
    \longmapsto
   (q,\eta+A(q)^{*}\lambda)
\end{align*}
where we identify elements in $\set{q}\times T_{[q]}^{*}(Q/K)$ with 
elements in $\hor^{*}_{q}$ via the dual of the inverse of the
horizontal lift. 

This \iso can be used to induce a symplectic form on the connection
dependent realization of $T^{*}Q$, namely $\sigma = \psi^{*}\Om$ where
$\Om=-d\theta$ is the canonical form on $T^{*}Q$. Moreover, there is
an explicit formula for $\sigma$ in terms of the chosen connection:
\[
 \sigma =
 (\wt{\pi}\circ\wt{\rho})^{*}\Om^{Q/K} - d{\dwt{\tau}}^{*}B
\]
where $\Om^{Q/K}$ is the canonical symplectic form on $T^{*}(Q/K)$, 
and furthermore
$B\in\Om^{1}(\bsc_{q}\ann\ko_{q})$ 
is given by 
\[
 B_{(q,\lambda)}(v_{1},\lambda_{1})
 = \vv<\lambda,A_{q}(v_{1})>.
\]
The explicit formula now is 
\begin{multline*}
 (dB)_{(q,\lam)}((v_{1},\lam_{1}),(v_{2},\lam_{2})) \\
  = 
 \vv<\lam,\curv^{A}_{q}(v_{1},v_{2})>
  + \vv<\lam,[Z_{1},Z_{2}]>
  - \vv<\lam_{2},Z_{1}> + \vv<\lam_{1},Z_{2}>
\end{multline*}
where
$(q,\lam)\in\bsc_{q}\ann\ko_{q}$, $(v_{i},\lam_{i})\in
T_{(q,\lam)}(\bsc_{q}\ann\ko_{q})$ for $i=1,2$, and 
\[
v_{i} 
  = 
  \zeta_{Z_{i}}(q)\oplus v_{i}^{\textup{hor}}
  \in
  \ver_{q}\oplus\hor_{q}
\]
is the decomposition into vertical and
horizontal part with $Z_{i}\in\ko$. 

Furthermore, there clearly is an induced action by $K$ on $\WW$.
This action is Hamiltonian with momentum
mapping 
\[
 \mu_{A} = \mu\circ\psi: \WW\longto\ko^{*}, (q,\eta,\lam)\longmapsto\lam,
\]
where $\mu$ is the
momentum map $T^{*}Q\to\ko^{*}$, 
and $\psi$ is equivariant.
\end{proposition}

\begin{proof}
This proposition is proved in Hochgerner \cite[Proposition
  5.1]{H04}. However, the proof can be considerably simplified using
the relation $\sigma = \psi^*\Om = -\psi^*d\theta$, and we 
present this simplification.

Let $\xi\in T\WW$ with $T(\wt{\tau}\circ\wt{\rho}).\xi = v\in TQ$ and $w=(q,\eta,\lam)\in\WW$.
Then 
\[
 (\psi^*\theta)_w(\xi) 
 =
 \vv<\eta+A_q^*(\lam),v_q>
 =
 ((\wt{\pi}\circ\wt{\rho})^{*}\theta^{Q/K})_w(\xi) 
 + 
 ({\dwt{\tau}}^{*}B)_w(\xi)
\]
whence
$
 \sigma 
 =
 (\wt{\pi}\circ\wt{\rho})^{*}\Om^{Q/K} - d{\dwt{\tau}}^{*}B
$.
It only remains to compute $dB$:
\begin{align*}
&(dB)_{(q,\lam)}((v_{1},\lam_{1}),(v_{2},\lam_{2}))\\
&=
\ld_{(v_{1},\lam_{1})} (B(v_2,\lam_2))_{(q,\lam)}
- \ld_{(v_{2},\lam_{2})} (B(v_1,\lam_1))_{(q,\lam)}
- B(q,\lam)([(v_{1},\lam_{1}),(v_{2},\lam_{2})])\\
&=
\ddo B(\fl_t^{(v_1,\lam_1)}(q,\lam))(v_2,\lam_2)
- \ddo B(\fl_t^{(v_2,\lam_2)}(q,\lam))(v_1,\lam_1)
- \langle \lam, A_q[v_1,v_2](q)\rangle\\
&=
\ddo \langle \fl_t^{\lam_1}(\lam),A(\fl_t^{v_1}(q))(v_2)\rangle
- \ddo \langle \fl_t^{\lam_2}(\lam),A(\fl_t^{v_2}(q))(v_1)\rangle
- \langle \lam, A_q[v_1,v_2](q)\rangle\\
&=
\langle \lam_1(\lam), Z_2\rangle 
- \langle \lam_2(\lam), Z_1\rangle
+ \langle \lam,\ld_{v_1}(A(v_2))_q-\ld_{v_2}(A(v_1))_q-A_q[v_1,v_2](q)\rangle\\
&= 
\langle \lam_1(\lam), Z_2\rangle 
- \langle \lam_2(\lam), Z_1\rangle
+ \langle \lam,d A_q(v_1,v_2)\rangle\\
&= 
\langle \lam_1(\lam), Z_2\rangle 
- \langle \lam_2(\lam), Z_1\rangle
+ \langle \lam,\curv^{A}_{q}(v_{1},v_{2})\rangle + \langle \lam,[Z_1,Z_2]\rangle.
\end{align*} 
\end{proof}

\section{Gauged Poisson reduction}\label{g-P-red}

Let $H$ be a closed subgroup of our compact Lie group $K$. Let $S$ be a
smooth manifold, and $a: H\times S\to S$, $(h,s)\mapsto h.s$ a left
action. We will be mostly interested in the case where $S$ is a slice
for the $K$-action on $Q$, and $a$ is the (trivial) slice
representation. 
Consider the left actions $l,r,t: H\times S\times K\to S\times K$ which are
defined by
\[
 l: (h,s,k)\longmapsto(s,hk),\text{ }
 r: (h,s,k)\longmapsto(s,kh^{-1}),\text{ and }
 t: (h,s,k)\longmapsto(h.s,kh^{-1}).
\]
The following lemma is elementary on the one hand as its computations are
straightforward. However, it is also tricky on the other hand since it
involves a choice of sign 
in the definition of the fundamental vectorfield from Section \ref{s:transfgr},
and there are many possibilities to get
confused. Note also that we use the left multiplication to trivialize
$T^{*}K=K\times\ko^{*}$. 

\begin{lemma}[Cotangent lifted actions]\label{lem:T^*K}
The lifted actions 
of $H$ on $T^{*}(S\times K)=T^{*}S\times K\times\ko^{*}$ 
corresponding to $l,r,t$ are given by 
\begin{align*}
  T^{*}S\times K\times\ko^{*}
 &\longto T^{*}S\times K\times\ko^{*}\\
T^{*}l_{h}: 
 (s,p,k,\eta)&\longmapsto(s,p,hk,\eta)\\
T^{*}r_{h}: 
 (s,p,k,\eta)&\longmapsto(s,p,kh^{-1},\Ad^{*}(h).\eta)\\
T^{*}t_{h}: 
 (s,p,k,\eta)&\longmapsto(h.s,T_{s}^{*}a_{h}.p,kh^{-1},\Ad^{*}(h).\eta)
\end{align*}
respectively, and where $\Ad^{*}(h).\eta :=
\eta\circ\Ad(h^{-1})$. Moreover, these lifted actions are Hamiltonian
with respect to the canonical exact symplectic form on 
$T^{*}(S\times K)=T^{*}S\times K\times\ko^{*}$, and the corresponding
momentum maps are given by
\begin{align*}
 T^{*}S\times K\times\ko^{*}&\longto\ho^{*}\\ 
 J_{l}: (s,p,k,\eta)&\longmapsto(\Ad^{*}(k).\eta)|\ho\\ 
 J_{r}: (s,p,k,\eta)&\longmapsto -\eta|\ho\\
 J_{t}: (s,p,k,\eta)&\longmapsto \mu(s,p)-\eta|\ho
\end{align*}
where $\mu: T^{*}S\to\ho^{*}$ is the canonical equivariant \momap with
respect to the Hamiltonian action $T^{*}a: H\times T^{*}S\to T^{*}S$. 
\end{lemma}

\begin{proof}
The point here is the choice of sign in the definition of
the fundamental vectorfield mapping $\zeta$ in Section \ref{s:transfgr}. Note
also that the cotangent bundle \momap $\mu: T^{*}S\to\ho^{*}$ is given
by $\vv<\mu(s,p),X>=\theta(\zeta_{X}^{a})(s,p)$ where $\theta$ is the
Liouville form on $T^{*}S$, and $\Om=-d\theta$ is the cotangent bundle
symplectic form.
\end{proof}

Let us introduce the abbreviation $E := \bsc_{q \in Q}\ann\ko_q$.

\begin{lemma}\label{E:Q}
The natural projection $\rho: E \to Q$ is a smooth  
fiber bundle with typical fiber $\ann\ho$, were $H$ is an isotropy 
subgroup of $K$ and $\ho$ is its Lie algebra.
\end{lemma}

\begin{proof}
We trivialize at an arbitrary point $q_0 \in Q$. 
We may assume $K_{q_0} = H$.
Let $U \subseteq Q$ be a tube around $K.q_0$ such that 
$U \cong S \times K/H$ as $K$-spaces, where $S$ is a slice at $q_0$.
Then it is possible to trivialize as
\[
\xymatrix{
 E|U  \ar[d] \ar[r]^-{\cong} & 
 S \times (K \times \ann\ho)/H \ar[d] \ar @{=}[r]  &
 S \times K \times_H \ann\ho  
             \\
 U \ar[r]^-{\cong} & S \times K/H &
}
\]
with trivializing map given by 
\[
 S\times K\times_H\ann\ho\longto E|U,
 \text{ }
 (s,[(k,\lambda)])\longmapsto(k.s,\Ad^*(k)(\lambda)) 
\]
where $\Ad^*(k)(\lambda) := \Ad(k^{-1})^*(\lambda)$. This map is 
well-defined and smooth with inverse given by
\[
(q,\lambda) = (k.s,\Ad^*(k)(\lambda_0)) \longmapsto (s,[(k,\lambda_0)])
\]
which is well-defined and smooth as well. Notice furthermore that the
trivializing map is $K$-equivariant with respect to the $K$-action on
$S\times K\times_{H}\ann\ho$ given by
$g.(s,[(k,\lam)])=(s,[(gk,\lam)])$. Indeed, this follows immediately
from Lemma \ref{lem:T^*K}.
\end{proof}

\begin{lemma}\label{lem:WW/K}
Let $E := \bsc_{q \in Q}\ann\ko_q$.
\begin{enumerate}[\up (i)]
\item
Let $U\subeq Q$ be a trivializing neighborhood for $\rho: E\to Q$ as
in the proof of Lemma \ref{E:Q}. Then $E|U$ is $K$-invariant, and if
$(L)$ is an element of the isotropy lattice for the $K$-action on $E$
then the corresponding stratum is trivialized as
\[
 (E|U)_{(L)}
 \cong
 S\times K\times_{H}(\ann\ho)_{(L_{0})^{H}}
\] 
where
$L_{0}\subeq H$ is an isotropy 
subgroup for the $H$-action conjugate to $L$ in $K$, 
and $(L_{0})^{H}$ is the conjugacy class of $L_{0}$ in $H$.
Moreover, the strata of $\ann\ho/H$ are of the form
$(\ann\ho)_{(L_{0})^{H}}/H$.
\item
The induced mapping $\rho_{0}: E/K\to Q/K$ is a singular fiber bundle
with typical fiber $\ann\ho/H$ in the sense of Definition
\ref{def:sg_fb}.
\end{enumerate}
\end{lemma}

\begin{proof}
Let $q_{0}\in Q$ with $K_{q_{0}}=H$. Then
\[
 (q_{0},\lam)\in(\bsc_{q\in Q}\ann\ko_{q})_{(L)} = E_{(L)}
\]
if and only if
\[
 \lam\in\ann\ho
 \textup{ and }
 H\cap K_{\lam}=H_{\lam}=L_{0}\sim L
 \textup{ within }K
\]
which is true if and only if
\[
 \lam\in(\ann\ho)_{(L_{0})^{H}}
\] 
where $L_{0}$ is a subgroup of $H$ conjugate to $L$ within $K$. 
Notice that it follows from the Slice Theorem for Riemannian actions
that $(\ann\ho)_{(L_{0})^{H}}$ is a smooth \mf, 
e.g. Palais and Terng \cite{PT88}.
Therefore, also the second assertion follows.
\end{proof}

\begin{theorem}\label{thm:WW/K}
There is a stratified isomorphism of stratified bundles over $Q/K$,
in the sense of Definition \ref{def:sg_fb},
\begin{align*}
 \psi_0^{-1} = \psi_0^{-1}(A): 
 (T^*Q)/K  &\longto  T^{*}(Q/K)\times_{Q/K}(\bsc_{q\in
                                Q}\ann\ko_{q})/K =: W,\\
 [(q,p)] &\longmapsto (C^{*}(q,p),[(q,\mu(q,p))])  
\end{align*}
where the stratification was suppressed. 
Here 
\[
 C^{*}: T^{*}Q\to\hor^{*}\to T^{*}(Q/K)
\] 
is constructed as 
the point wise dual to the horizontal lift mapping $C:
T(Q/K)\times_{Q/K}Q\to\hor\subeq TQ$, $([q],v;q)\to C_{q}(v)$ 
associated to the connection $A \in \Omega^1(Q;\ko)$. 

If $(L)$ is an isotropy class of the $K$-action on $T^{*}Q$, then 
$\psi_{0}^{-1}$ maps the isotropy stratum $(T^*Q)_{(L)}/K$ onto 
\[
 T^{*}(Q/K)\times_{Q/K}(\bsc_{q \in Q}\ann\ko_{q})_{(L)}/K =: W_{(L)}.
\]
Moreover, the natural projection 
\[
 \widetilde{\rho}_{0}^{(L)}: W_{(L)} \to T^*(Q/K)
\]
is a smooth fiber bundle with typical fiber of the form 
$(\ann\ho)_{(L_{0})^{H}}/H$.
Here $L_{0}\subeq H$ is an isotropy 
subgroup for the $H$-action conjugate to $L$ in $K$, 
and $(L_{0})^{H}$ is the conjugacy class of $L_{0}$ in $H$.

Therefore, $\widetilde{\rho}_{0}: W\to T^{*}(Q/K)$ is a singular 
fiber bundle in the sense of Definition \ref{def:sg_fb}. 
\end{theorem}

In the case that $K$ acts on $Q$ freely the first assertion of the
above theorem can also be
found in Cendra, Holm, Marsden, Ratiu \cite{CHMR98}. Following Ortega
and Ratiu \cite[Section 6.6.12]{OR04} the above constructed
interpretation $W$ of $(T^{*}Q)/K$ is called \caps{Weinstein space}
referring to Weinstein \cite{Wei78} where this universal
construction first appeared. 

\begin{proof}
We consider first the map 
\begin{align*}
 \varphi_0 = \varphi_0(A): 
 (TQ)/K  &\longto  T(Q/K)\times_{Q/K}(\bsc_{q\in
                                Q} \ko_{q}^{\bot})/K,\\
 [(q,v)] &\longmapsto (T\pi(q,v),[(q,A_q v)])  
\end{align*}
the point-wise dual of whose inverse will be $\psi_0^{-1}$. 

The spaces $TQ$ and $(TQ)/K$ are stratified into isotropy types. 
Since the base $Q$ is stratified as consisting only of a single stratum, 
the equivariant foot point projection map $\tau: TQ\to Q$ is 
trivially a stratified map. 
Using the Slice Theorem on the base $Q$ it is easy to see that both
$(TQ)/K\to Q/K$ 
and the projection 
$(\bsc_{q\in Q}\ko_{q}^{\bot})/K\to Q/K$ are singular bundle
maps in the sense of Definition \ref{def:sg_fb}.
Hereby
$(\bsc_{q\in Q}\ko_{q}^{\bot})/K$ is stratified into
isotropy types.
According to Davis \cite{Dav78} or also Subsection \ref{sec:pullbacks}
pullbacks are well defined in
the category of stratified spaces and thus it makes sense to define
$T(Q/K)\times_{Q/K}(\bsc_{q\in Q}\ko_{q}^{\bot})/K$ as a stratified
space with smooth structure.

The map $\varphi_0$ is well defined: indeed, for $(q,v)\in TQ$ and 
$k\in K$ 
we have 
\[
 T\pi(k.q,k.v) 
  = (\pi(k.q),T_{k.q}\pi(T_{q}l_{k}(v))) 
  = (\pi(q),T_{q}(\pi\circ l_{k})(v)) 
  = T_{q}\pi(v),
\] 
and
$[(k.q,A(k.q,k.v))] = [(q,A(q,v))]$ 
by equivariance of $A$.
It is clearly continuous as a composition of continuous maps. 
Moreover, since $\cinf((TQ)/K)=\cinf(TQ)^{K}$ by Example
\ref{ex:quot_smooth} it follows that $\varphi_0$ is a smooth map of
singular spaces. 

We claim that $\varphi_0$ maps strata onto strata, and moreover we
have the formula
\[
 \varphi_0((TQ)_{(L)}/K)
 =
 T(Q/K)\times_{Q/K}(\bsc_{q\in Q}\ko_{q}^{\bot})_{(L)}/K.
\]
Indeed, consider $(q,v)\in(TQ)_{(L)}$, that is 
$H\cap K_{v} = L^{\prime}\sim L$ where $H=K_{q}$. 
The notation $L^{\prime}\sim L$ means that $L^{\prime}$ is conjugate 
to $L$ within $K$.
Now we can decompose $v$ as
$v = v_{0}\oplus\zeta_{X}(q)\in\hor_{q}\oplus T_{q}(K.q)$ for some
appropriate $X\in\ko$.
Since $Q$
consists only of a single isotropy type we have $T_{q}Q =
T_{q}Q_{H}+T_{q}(K.q)$ -- which is not a direct sum
decomposition. As usual, $Q_{H}=\set{q\in Q: K_{q}=H}$.
This shows that $v_{0}\in T_{q}Q_{H}$, 
since $\hor_q = T_q(K.q)^{\bot} \subseteq T_q Q_H$, 
and hence $H\subeq K_{v_{0}}$. 
By equivariance of $A$ it follows that 
\[
 K_{q}\cap K_{v}
 = 
 H\cap K_{v_{0}}\cap K_{\zeta_{X}(q)}
 =
 H\cap K_{\zeta_{X}(q)}
 =
 H\cap K_{A(q,v)}
\]
which is independent of the horizontal component. Hence the
claim. The restriction of $\varphi_0$ to any stratum clearly is smooth
as a composition of smooth maps. 

Since $A_{q}(\zeta_{X}(q))=X$ for
$X\in\ko_{q}^{\bot}$ we can write down an inverse as
\[
 \varphi_0^{-1}:
 ([q],v;[(q,X)])
 \to
 [(q,C_{q}(v)+\zeta_{X}(q))]
\]
and again it is an easy matter to notice that this map is well
defined, continuous, and smooth on each stratum. 
Again, it follows from the definition of the smooth structures on the
respective spaces that $\varphi_0^{-1}$ is smooth.

It makes sense to define the dual $\psi_0^{-1}$ of the inverse map
$\varphi_0^{-1}$ in a point wise manner, and it only remains to compute
this map.
\begin{align*} 
 \vv<\psi_0^{-1}[(q,p)],([q],v;[(q,X)])> 
 &= 
 \vv<[(q,p)],[(q,C_{q}(v)+\zeta_{X}(q))]> \\
 &= 
 \vv<p,C_{q}(v)> + \vv<p,\zeta_{X}(q)> \\
 &=
 \vv<C^{*}(q,p),v> + \vv<\mu(q,p),X> \\
 &=
 \vv<(C^{*}(q,p),[(q,\mu(q,p))]),([q],v;[(q,X)])>
\end{align*}
where we used the $K$-invariance of the dual pairing over $Q$. 

Finally, $\psi_0^{-1}$ is an \iso of singular Poisson spaces: 
note first that the identifying map 
$\WW/K\to W$, 
$[(q;[q],\eta;q,\lam)]\mapsto ([q],\eta;[(q,\lam)])$ 
is well-defined because $K_{q}$ acts
trivially on $\hor_{q}^{*}=T^{*}_{[q]}(Q/K)\ni\eta$ which in turn is due
to the fact that all points of $Q$ are regular. 
Moreover, by the universal property for singular pull back bundles
from Subsection \ref{sec:pullbacks} it is obvious that this map 
$\WW/K\to W$ is smooth and has a smooth inverse.
The quotient Poisson
bracket is well-defined since $\cinf(\WW)^{K}\subeq\cinf(\WW)$ is a
Poisson sub-algebra. The statement now follows because the diagram
\[
\xymatrix{
 T^{*}Q \ar[r]^-{\psi^{-1}}\ar @{->>}[d] & {\WW}\ar @{->>}[dr] \\
 (T^{*}Q)/K\ar[r]^-{\psi_0^{-1}} & W\ar @{=}[r] & {\WW}/K
}
\]
is commutative, and composition of top and down-right arrow is
Poisson and the left vertical arrow is surjective. 

Using Lemma \ref{lem:WW/K} together with
Subsection \ref{sec:sg_fb} the global description of $\WW/K\cong W$ 
as a fibered product thus follows. 
\end{proof}

Next we shall construct a connection on 
$\rho : E \to Q$ which will provide a connection on 
$\widetilde{\rho} : \WW \to Q \times_{Q/K} T^*(Q/K)$. 
Recall the mechanical connection $A\in\Om^{1}(Q;\ko)$ from
Section \ref{W-bundle}.
Consider the embedding
\[
 \iota : E\longto T^*Q,\text{ }
 (q,\lambda)\longmapsto (q,A_q^*(\lambda)).
\]
On $\varrho : T^*Q \to Q$ we choose the canonical (with respect to 
the metric) linear connection $\Phi(\varrho) : TT^*Q \to V(\varrho)$. 
Consider the following diagram 
\begin{align} \label{diag}
\xymatrix{
 TE \ar@{-->}[r]^-{\Phi(\rho)} \ar@{^{(}->}[d]_{T \iota} &
  {V(\rho)}\ar @{=}[r]
  &{\bsc_{(q,\lam)\in E}\ann\ko_{q}}\\
 TT^*Q \ar[r]^{\Phi(\varrho)} &
 V(\varrho) \ar[u]_{d \mu |{V(\varrho)}}
}
\end{align}
which induces a linear connection $\Phi(\rho)$ on 
$\rho: E \to Q$. 

\begin{lemma}
The Diagram \eqref{diag} is well-defined.
\end{lemma}

\begin{proof}
We have to show that $d \mu |{V(\varrho)}$ takes values in 
$V(\rho) = \bsc_{(q,\lam)\in E}\ann\ko_{q}$. 
Since $\mu_q(A_q^*(\lambda)) = \lambda$ for all $\lambda \in
\ann\ko_q$ and $q\in Q$, 
it suffices to prove that $\vv<d\mu(q_0,p_0)(\xi_{(q_0,p_0)}),X> = 0$
for fixed
$(q_0,p_0)\in T^*Q$, $\xi \in \Gamma(V(\varrho))$, and $X \in\ko_{q_0}$. 
Now, for arbitrary $(q,p)\in T^*Q$, 
\begin{eqnarray*}
 \vv<d\mu(q,p)(\xi_{(q,p)}),X>
  &=& 
 \Omega_{(q,p)}(\zeta_X^{T^*Q},\xi) \\ 
  &=& 
  - (\zeta_X^{T^*Q}.\theta(\xi))(q,p) 
  + (\xi.\theta(\zeta_X^{T^*Q}))(q,p)
  + \theta([\zeta_X^{T^*Q},\xi])(q,p).
\end{eqnarray*}
We have 
\[
\theta(\xi)(q,p) = \langle p , T_{(q,p)} \varrho.\xi_{(q,p)} \rangle = 0
\] 
for all $(q,p) \in T^*Q$ and $\xi \in \Gamma(V(\varrho))$, so $\theta(\xi) = 0$. 
Next we find 
\[
\theta(\zeta_X^{T^*Q})(q_0,p) 
 = \vv<p,T_{(q_0,p)}\varrho.\zeta_X^{T^*Q}(q_0,p)>
 = \vv<p,\zeta_X(q_0)>
 = 0
\]
for all $p\in T_{q_0}^*Q$ since $X \in \ko_{q_0}$ by assumption.
Since $\xi$ takes values in the vertical subbundle $V(\varrho)$ and
vertical bundles are integrable its flow preserves fibers of
$\varrho$, that is $\fl^{\xi}_t(q_0,p_0)=(q_0,p_t)$ for some curve
$p_t$ in $T_{q_0}^*Q$. Therefore,
\[
 (\xi.\theta(\zeta_X^{T^*Q}))(q_0,p_0)
 = \dd{t}{}|_0\vv<p_t,\zeta_X(q_0)>
 = \dd{t}{}|_0 0 
 = 0.   
\]
Finally, 
\[
\theta([\zeta_X^{T^*Q},\xi])(q,p) 
= \langle p , T_{(q,p)} \varrho.[\zeta_X^{T^*Q},\xi](q,p) \rangle 
= 0
\]
for all $(q,p) \in T^*Q$, since the fact that $\zeta_X^{T^*Q}$ and 
$\zeta_X$ are $\varrho$-related as well as $\xi$ and $0$ implies that 
$[\zeta_X^{T^*Q},\xi]$ and $[\zeta_X,0]=0$ are $\varrho$-related, 
i.e.\ $T_{(q,p)} \varrho.[\zeta_X^{T^*Q},\xi](q,p) = 0$. 
This completes the proof.
\end{proof}

Via the pullback 
construction (see Diagram \eqref{W-diag}) 
this also induces a linear connection on 
$\widetilde{\rho} : \WW \to Q \times_{Q/K} T^*(Q/K)$. We denote 
this connection by $\widetilde{A} : T\WW \to V(\widetilde{\rho})$.
Notice that $V_{(q,\eta,\lam)}(\widetilde{\rho}) = \ann\ko_{q}$ 
by construction.

The connection $\widetilde{A}$ and the momentum map 
$\mu_A := \mu \circ \psi : \WW \to \ko^*$ are related by
\[
\widetilde{A}(q,\eta,\lambda)(\xi) = d \mu_A(q,\eta,\lambda)(\xi),
\] 
where $\xi \in T_{(q,\eta,\lambda)}\WW$, and $(q,\eta,\lambda)$ is 
short-hand for $(q;[q],\eta;q,\lambda) \in \WW$.

We will use the connection $\widetilde{A}$ to decompose an arbitrary 
vector $\xi \in T_{(q,\eta,\lambda)} \WW$ as
\[
\xi = (v(q);\eta'([q],\eta);v_1(q),\nu(q,\lambda)), 
\]
where $\nu(q,\lambda) = \widetilde{A}(q,\eta,\lambda)(\xi)$ is 
independent of $\eta$. 
Notice also that $v_1(q) = v(q)$ by the pullback property. 
Further we can decompose $v(q) \in T_q Q$ according to
\[
v(q) = v(q)^{\textup{hor}(\pi)} + \zeta_Z(q) \in H_q(\pi) \oplus V_q(\pi)
\]
with respect to the connection $A$ on $\pi : Q \to Q/K$.
The same can be done with 
$\eta'([q],\eta) \in T_{([q],\eta)}(T^*(Q/K))$ as
\[
\eta'([q],\eta)  = \eta'([q],\eta)^{\textup{hor}(\tau)} 
                 + \eta'([q],\eta)^{\textup{ver}(\tau)}
                 \in H_{([q],\eta)}(\tau) \oplus V_{([q],\eta)}(\tau)
\]
with respect to the canonical connection on $\tau : T^*(Q/K) \to Q/K$ 
which comes from the induced metric on $Q/K$. Notice that we have 
$\eta'([q],\eta)^{\textup{hor}(\tau)} = v(q)^{\textup{hor}(\pi)}$ 
by the pullback property.

\begin{definition}[Vertical differentiation on 
$\widetilde{\rho} : \WW \to Q \times_{Q/K} T^*(Q/K)$]
Let $\nu \in V_{(q,\eta,\lambda)}(\widetilde{\rho}) = \ann\ko_{q}$ and 
$F \in C^{\infty}(\WW)$.
We define 
\[
d_v F(q,\eta,\lambda)(\nu) := \left.\frac{\partial}{\partial t}\right|_0
                              F(q,\eta,\lambda+t \nu)
\]
to be the \caps{vertical derivative} of $F$ at $(q,\eta,\lambda)$.
\end{definition}

\begin{definition}[Covariant differentiation on 
$\widetilde{\rho} : \WW \to Q \times_{Q/K} T^*(Q/K)$]
Let $\chi := 1-\widetilde{A} : T \WW \to H(\widetilde{\rho})$ 
denote the horizontal projection with respect to $\widetilde{\rho}$. 
The \caps{covariant derivative} of $F \in C^{\infty}(\WW)$ is defined 
as $d_{\widetilde{A}}F := d F \circ \chi$.
\end{definition}

\begin{lemma} \label{lem1}
Let $F \in C^{\infty}(\WW)^K$, and decompose the Hamiltonian vector field of $F$ at 
$(q,\eta,\lambda) \in \WW$ as
\[
\nabla_F^{\sigma}(q,\eta,\lambda) 
= (v(q);\eta'([q],\eta);v(q),\nu(q,\lambda))
\]
with 
$\nu(q,\lambda) 
= \widetilde{A}(q,\eta,\lambda)(\nabla_F^{\sigma}(q,\eta,\lambda))$
according to above. Here $\sigma$ is the symplectic structure on 
$\WW$ from Proposition \ref{prop:WW}. Then, $\nu(q,\lambda) = 0$.
\end{lemma}

\proof
This is a consequence of Noether's Theorem. We have
\[
\nu(q,\lambda) 
= \widetilde{A}(q,\eta,\lambda)(\nabla_F^{\sigma}(q,\eta,\lambda))
= d \mu_A(q,\eta,\lambda)(\nabla_F^{\sigma}(q,\eta,\lambda))
= 0,
\]
since $\mu_A$ is constant along flow lines of Hamiltonian vector fields 
of invariant functions.
\endproof

\begin{lemma} \label{lem2}
Let $F \in C^{\infty}(\WW)^K$, and decompose the Hamiltonian vector field of 
$F$ at $(q,\eta,\lambda) \in \WW$ as
\[
\nabla_F^{\sigma}(q,\eta,\lambda) 
= (v(q);\eta'([q],\eta);v(q),0).
\]
Then,
\[
\eta'([q],\eta) 
= \left.(\Omega_{([q],\eta)}^{Q/K} \check{)} \right.^{-1}
(d_{\widetilde{A}} F(q,\eta,\lambda)), 
\]
where we consider $d_{\widetilde{A}} F(q,\eta,\lambda)$ as an element 
of $T^*_{([q],\eta)}(T^*(Q/K))$ via the dual of the horizontal lift mapping. 
\end{lemma}

Recall that $\Omega^{Q/K}$ denotes 
the canonical symplectic form on $T^*(Q/K)$.

\begin{proof}
For the purpose of the lemma we may assume that the vertical part of
the first entry vanishes, i.e\
$A.T(\widetilde{\tau} \circ \widetilde{\rho}).\nabla^{\sigma}_F = 0$ 
with notation as in Diagram \eqref{W-diag}.
Thus we may simplify the problem by assuming that
$F\in\cinf(\WW)^K$ factors through 
$\widetilde{\pi} \circ \widetilde{\rho}$, i.e.\ there is a smooth map
$f: T^{*}(Q/K)\to\R$ such that 
$F =
f\circ\widetilde{\pi}\circ\widetilde{\rho}: \WW\to\R$.
Then, 
$d_{\widetilde{A}} F(q,\eta,\lam) 
= d_{\widetilde{A}} F(q,\eta,0) 
= d f((\widetilde{\pi} \circ \widetilde{\rho})(q,\eta,0)) 
= d f([q],\eta)$. 
Now, the projection 
\[
\xymatrix{
 {
 Q\times_{Q/K}T^*(Q/K)
 =
 \mu_{A}^{-1}(0)}
 \ar @{->>}[r]
 & 
{\mu_A^{-1}(0)/K = T^*(Q/K)}
}
\]
is a Poisson morphism. 
Therefore, we find
\begin{align*}
\eta'([q],\eta) 
 &= T_{(q,\eta,0)}(\widetilde{\pi} \circ \widetilde{\rho})
      .\nabla_F^{\sigma}(q,\eta,0)
  = \nabla_f^{Q/K}((\widetilde{\pi} \circ \widetilde{\rho})
                 (q,\eta,0)) \\
 &= \left.(\Omega_{([q],\eta)}^{Q/K} \check{)} \right.^{-1}
        (d f([q],\eta))
  = \left.(\Omega_{([q],\eta)}^{Q/K} \check{)} \right.^{-1}
      (d_{\widetilde{A}} F(q,\eta,0))
\end{align*}
which gives the assertion.
\end{proof}

\begin{lemma} \label{lem3}
Let $F \in C^{\infty}(\WW)^K$, and decompose the Hamiltonian vector field of 
$F$ at $(q,\eta,\lambda) \in \WW$ as
\[
\nabla_F^{\sigma}(q,\eta,\lambda) 
= (v(q);\eta'([q],\eta);v(q),0).
\]
Via the connection $A \in \Omega^1(Q;\ko)$ we can write 
\[
v(q) = v(q)^{\textup{hor}(\pi)} + \zeta_Z(q) \in H_q(\pi) \oplus V_q(\pi).
\]
Then $Z = -d_v F(q,\eta,\lam)\in\ko_{q}^{\bot}$, and, moreover,
$\Ad(h).Z = Z$ for all $h\in K_{q}$.
\end{lemma}

\begin{proof}
We work in tube coordinates around $q \in Q$. Thus let $S$ be a slice 
through $q$ for the $K$-action on $Q$ such that $U \cong S \times K/H$
as $K$-spaces 
where $U$ is a tube around $K.q$ and $H=K_q$. 
Then we have a $K$-equivariant isomorphism of symplectic manifolds 
\[
\WW|U \cong T^{*}S \times K \times_H \ann\ho
 \cong T^{*}S\times T^{*}K\spr{0}T^*R(H)
\]
where the right hand side carries the obvious symplectic structure.
Here $T^*R(H)$ is the cotangent lifted action of the right multiplication 
of $H$ on $K$, and we use left multiplication to trivialize 
$T^*K = K\times\ko^*$.
This follows from Lemma \ref{lem:T^*K} and an argument similar as in
the proof of Lemma \ref{E:Q}. In particular, note that (as in
Hochgerner \cite[Section 4]{H04}) the isomorphism is symplectic since
it comes as the cotangent lift of a diffeomorphism of the base spaces.
Since we already know that the part of the Hamiltonian vector field of $F$ 
which is tangent to $T^*S$
is given by the local coordinates of
$(v(q)^{\textup{hor}(\pi)},\eta'([q],\eta))$, 
we may further reduce the problem to considering 
a function 
$F \in C^{\infty}(K \times_H \ann\ho)^K = C^{\infty}(\ann\ho/H)$. 
This identification is due to Lemma \ref{lem:WW/K}(ii).
Now, referring again to Lemma \ref{lem:T^*K} we have
\[
 K\times_H\ann\ho = T^*K\spr{0}T^*R(H),
\]
where $T^*R(H)$ is the cotangent lifted action of the right multiplication 
of $H$ on $K$, and where we use left multiplication to trivialize 
$T^*K = K\times\ko^*$.
Thus, there exists a function 
$F_1 \in C^{\infty}(T^*K) = C^{\infty}(K \times \ko^*)$
that is $T^*L(K)$-invariant ($L$ denotes the left multiplication on $K$) 
and $T^*R(H)$-invariant
such that the following diagram commutes:
\[
\xymatrix{
K \times \ann\ho \ar@{^{(}->}[r] \ar@{->>}[d] &
K \times \ko^* \ar[d]^{F_1} \\
K \times_H \ann\ho \ar[r]^-{F} &
{\R}
}
\]
We choose local cotangent bundle coordinates 
$a_i,b_i$ where $i = 1,\ldots,m$ on 
$T^*K = K \times \ko^*$
such that 
$b_1,\ldots,b_l$ 
are coordinates on $\ho^*$,
$b_{l+1},\ldots,b_m$ 
are coordinates on $\ann\ho$,
and such that 
$\by{\partial}{\partial a_1},\ldots,\by{\partial}{\partial a_l}$ 
are a basis of $\ho$, and
$\by{\partial}{\partial a_{l+1}},\ldots,\by{\partial}{\partial a_m}$ 
are a basis of $\ho^{\bot}$. 
Then for the canonical Poisson bracket on $T^*K$ we obtain
\begin{multline*}
 -\nabla^{T^{*}K}_{F_{1}}(e,\lam)
 =
 \{F_1,\cdot\}^{T^*K}(e,\lambda) 
 =\\
 = 
 \sum_{i=1}^m 
 \left( 
  \by{\partial F_1}{\partial b_i}(e,\lambda)\by{\partial}{\partial a_i} 
  -
  \by{\partial F_1}{\partial a_i}(e,\lambda)\by{\partial}{\partial b_i}
 \right)
 =
 \sum_{i=l+1}^m 
 \by{\partial F_1}{\partial b_i}(e,\lambda)\by{\partial}{\partial a_i} 
 \in\ho^{\bot}
\end{multline*}
which is the vertical derivative of $F_1$ identified with an element 
of $\ho^{\bot}$ through choice of a basis. 
Since the projection $K \times \ann\ho \toto K \times_H \ann\ho$ 
is Poisson, Hamiltonian vector fields project to Hamiltonian vector 
fields,
i.e.\ $-\nabla^{T^{*}K}_{F_{1}}(e,\lam)$ projects to 
\[
 -\nabla^{K\times_{H}\ann\ho}_{F}[(e,\lam)] = -Z.
\]
Therefore, the Hamiltonian vector field of $F$ on 
$K\times_H\ann\ho$ is $(d_v F,0)$ which is tangent to the $K/H$-factor.
Thus $-Z = d_v F(q,\eta,\lam)\in\ko_q^{\bot}$.

To see the second assertion let $h\in H$ and 
notice that 
\[
 \Ad(h).\nabla^{T^{*}K}_{F_{1}}(e,\lam)
 =
 T_{h}R_{h^{-1}}.T_{e}L_{h}.\nabla_{F_{1}}^{T^{*}K}(e,\lam)
 =
 \nabla_{F_{1}}^{T^{*}K}(e,\Ad^*(h).\lam)
\]
since $F_1$ is both $T^*L(H)$- and $T^*R(H)$-invariant. 
Let $\pi_H: K\times\ann\ho\toto K\times_H\ann\ho$ denote the
projection. Then
\[
 \Ad(h).Z
 =
 T\pi_H.\nabla_{F_{1}}^{T^{*}K}(e,\Ad^*(h).\lam)
 =
 \nabla^{K\times_{H}\ann\ho}_{F}[(h,\lam)]
 =
 T\pi_H.\nabla_{F_{1}}^{T^{*}K}(e,\lam)
 = 
 Z
\]
since $F_1$ depends only on the second factor, i.e.\ it is
$T^*L(K)$-invariant.
\end{proof}

\begin{theorem}[Poisson structure on Weinstein space]\label{thm:p-struct}
The identification
\begin{align*}
 \WW/K \overset{=}{\longto} W, \text{ }
 [(q;[q],\eta;q,\lam)] \longmapsto ([q],\eta;[(q,\lam)])
\end{align*}
gives an induced Poisson bracket on $C^{\infty}(W) = C^{\infty}(\WW)^K$ 
which makes the stratified isomorphism
\[
\psi_0 = \psi_0(A): 
 T^{*}(Q/K)\times_{Q/K}(\bsc_{q\in
                                Q}\ann\ko_{q})/K = W \longto (T^*Q)/K
\]
from Theorem \ref{thm:WW/K} into an isomorphism of Poisson spaces.

Let $[(q,\eta,\lambda)] \in \WW/K = W$, and $f_1,f_2 \in C^{\infty}(W)$. 
Assume that $F_1,F_2 \in C^{\infty}(\WW)^K$ are lifts of $f_1,f_2$ to 
$\WW$. Then the induced Poisson bracket on $W$ is given by
\begin{align*}
 &\bra{f_1,f_2}^W[(q,\eta,\lam)]\\ 
 &\phantom{xx}=  
   \Omega_{([q],\eta)}^{Q/K}
    \left( 
    \left.(\Omega_{([q],\eta)}^{Q/K} \check{)} \right.^{-1}
           (d_{\widetilde{A}} F_1(q,\eta,\lambda)),
    \left.(\Omega_{([q],\eta)}^{Q/K} \check{)} \right.^{-1}
           (d_{\widetilde{A}} F_2(q,\eta,\lambda))
    \right)\\
 &\phantom{xx=}- 
  \vv<\lam,\curv_{0}^{A}(v_1(q)^{\textup{hor}(\pi)},v_2(q)^{\textup{hor}(\pi)})>
  - 
  \vv<\lam,[d_v F_1(q,\eta,\lam),d_v F_2(q,\eta,\lam)]>
\end{align*}
where $\Omega^{Q/K}$ is the canonical symplectic form on $T^*(Q/K)$, 
the covariant derivatives $d_{\widetilde{A}} F_i(q,\eta,\lambda)$ 
are considered as elements 
of $T^*_{([q],\eta)}(T^*(Q/K))$, 
and 
$v_i(q)^{\textup{hor}(\pi)} 
      = T\tau.\left.(\Omega_{([q],\eta)}^{Q/K} \check{)}
         \right.^{-1}(d_{\widetilde{A}}F_i(q,\eta,\lambda))$. 
Finally, 
$\curv_{0}^{A}$ is the 
induced form on $Q/K$ associated to the mechanical connection $A$ from
Proposition \ref{prop:curvA}.
\end{theorem}

In the case that $K$ acts freely on $Q$ the Poisson bracket on the 
reduced Poisson manifold $T^*Q/K$ is determined in Zaalani \cite{Zaa99} 
and in Perlmutter and Ratiu \cite{PR04}. In the first paper the realization 
of $T^*Q/K$ as Weinstein space is used, the latter deals with its 
realization as Sternberg and Weinstein space. 

\begin{proof}
The first part of the theorem has already been checked in the proof of 
Theorem \ref{thm:WW/K}.

Let $f_1,f_2 \in C^{\infty}(W)$ and let $F_1,F_2 \in C^{\infty}(\WW)^K$ be its unique
lifts to $\WW$. 
In order to establish the formula for the reduced Poisson bracket we decompose
the Hamiltonian vector fields of $F_1$ and $F_2$ as above
\[
\nabla_{F_i}^{\sigma}(q,\eta,\lambda) 
= (v_i(q);\eta_i'([q],\eta);v_i(q),\nu_i(q,\lambda))
\qquad (i=1,2).
\]
With the intrinsic symplectic form $\sigma$ on $\WW$ from Proposition \ref{prop:WW} 
we have
\[
 \{f_1,f_2\}^W [(q,\eta,\lam)] 
   = \{F_1,F_2\}^{\mathcal{W}} (q,\eta,\lam) 
   = \sigma(\nabla_{F_1}^{\sigma},\nabla_{F_2}^{\sigma})(q,\eta,\lam)
\]
which turns to the desired formula by the identity 
$\eta_i'([q],\eta)^{\textup{hor}(\tau)} = v_i(q)^{\textup{hor}(\pi)}$, 
and Lemmas \ref{lem1}, \ref{lem2}, and \ref{lem3}.
\end{proof}

The above theorem implies, in particular, that the embedding
$T^*(Q/K)\hookto W$ as the zero section is a Poisson morphism when
$T^*(Q/K)$ is equipped with its standard Poisson structure. Now,
fixing $\lam\in(\ann\ho)_{(L)^H}$ defines a smooth embedding
\[
 \iota_{\lam}:
 T^*(Q/K)\cong T^*(Q_H/N(H))
 \longto
 W_{(L)},
 \text{ }
 ([q],\eta)\longmapsto[(q,\eta,\lam)]
\]
where $q\in Q_H := \set{q\in Q: K_q=H}$, $N(H)$ is the normalizer of $H$ in $K$, 
and $(L)^H$, $(L)$ denotes
the conjugacy class of $L\subeq H$ in $H$, $K$ respectively.  
A quick inspection shows that $\iota_{\lam}$ is, in general,
not Poisson if
$T^*(Q/K)$ carries its usual Poisson structure. An element
$\lam\in\ko^*$ is called totally isotropic if $\ko_{\lam}=\ko$.

\begin{corollary}[Charge]\label{cor:charge}
If $\lam\in\ann\ho$ is totally isotropic then the embedding
$\iota_{\lam}: T^*(Q/K)\hookto W$ is Poisson if $T^*(Q/K)$ is equipped
with the Poisson structure stemming from the magnetic symplectic form
$\Om^{\lam} := \Om^{Q/K}-\vv<\lam,\tau^*\curv^{A}_{0}>$. 
\end{corollary}

\begin{proof}
This is obvious from the formula of Theorem \ref{thm:p-struct}.
\end{proof}

Due to the appealing similarity of the symplectic form
$\Om^{\lam}$ with that appearing in electromagnetism $\Om^{\lam}$ is
called a magnetic symplectic form, and we think of the totally
isotropic momentum value $\lam$ as the charge of the 
testparticles moving in the electromagnetic field on the
reduced configuration space $Q/K$. 
This
analogy will be carried further in Section \ref{symp-leaves} where we
equip our particles with spin by considering general momentum
values. See Guillemin and Sternberg \cite[Section 20]{GS84} for a
discussion of electromagnetism in a symplectic framework.

Since Hamiltonian vector fields associated to $K$-invariant functions
are tangent to the isotropy type submanifolds  we can define
Hamiltonian vector fields on $W$ as follows. Let $f\in\cinf(W)$ and
$F=\phi^*f\in\cinf(\WW)^K$ where $\phi: \WW\toto W$ is the
projection. Then $\nabla^W_f$ which is characterized by 
\[
 \nabla^W_f \circ \phi 
 = 
 -\bra{f,\cdot}^W \circ \phi
 = 
 T\phi.\nabla^{\sigma}_F
\]
is a stratified \vf on $W$ in the sense of Subsection
\ref{sec:sg_fb}. Let $\phi_{(L)}: \WW_{(L)}\toto W_{(L)}$ denote the
restriction of $\phi$ to the isotropy stratum. 

\begin{corollary}[Hamiltonian vector fields]\label{cor:hvf}
Let $[(q,\eta,\lam)]\in W_{(L)} = \WW_{(L)}/K$ and
$f\in\cinf(W)$ with unique lift $F=\phi^*f$. 
Then
\[
 \nabla_f^W[(q,\eta,\lam)]
 =
 (v_0,
  \left.(\Omega_{([q],\eta)}^{Q/K} \check{)}
         \right.^{-1}(d_{\widetilde{A}}F(q,\eta,\lambda)),
  \ad^*(d_{v}F(q,\eta,\lam)).\lam)
\]
where we consider $d_{\widetilde{A}}F(q,\eta,\lambda)$ as an element
of $T^*_{([q],\eta)}(T^*(Q/K))$ through the isomorphism given by the
dual of the horizontal lift with respect to the mechanical connection
$A$ on $Q\toto Q/K$. Moreover, 
$v_0 := 
 T_{([q],\eta)}\tau.\left.(\Omega_{([q],\eta)}^{Q/K} \check{)}
         \right.^{-1}(d_{\widetilde{A}}F(q,\eta,\lambda))$.
\end{corollary}


\begin{proof}
Notice firstly that 
\[
 \nabla^{\sigma}_F(q,\eta,\lam)
 -
 \Big(v(q)^{\textup{hor}(\pi)},
  \left.(\Omega_{([q],\eta)}^{Q/K} \check{)}
         \right.^{-1}(d_{\widetilde{A}}F(q,\eta,\lambda)),
  \ad^*(-Z).\lam\Big)
 \in \ker T_{(q,\eta,\lam)}\phi_{(L)}
\]
by Lemmas \ref{lem1}, \ref{lem2}, and \ref{lem3}, 
and where $-Z = d_{v}F(q,\eta,\lam)$.
We have $v_0 =
T_{q}\pi.(v(q)^{\textup{hor}(\pi)})$
then
\begin{align*}
 \nabla_f^W[(q,\eta,\lam)]
 &=
  T_{(q,\eta,\lam)}\phi_{(L)}.\nabla^{\sigma}_F(q,\eta,\lam)\\
 &=
  T_{(q,\eta,\lam)}\phi_{(L)}.
   (v(q)^{\textup{hor}(\pi)},
   \left.(\Omega_{([q],\eta)}^{Q/K} \check{)}
          \right.^{-1}(d_{\widetilde{A}}F(q,\eta,\lambda)),
   \ad^*(-Z).\lam)\\
 &=
  (v_0,
   \left.(\Omega_{([q],\eta)}^{Q/K} \check{)}
          \right.^{-1}(d_{\widetilde{A}}F(q,\eta,\lambda)),
    \ad^*(-Z).\lam).
\end{align*}
The last equality in the aligned
equation is true since 
\[
 (v(q)^{\textup{hor}(\pi)},
   \left.(\Omega_{([q],\eta)}^{Q/K} \check{)}
          \right.^{-1}(d_{\widetilde{A}}F(q,\eta,\lambda)),
   \ad^*(-Z).\lam)
\]
is horizontal with respect to the projection
$\phi_{(L)}: \WW_{(L)}\toto W_{(L)}$ which is equipped with the structure
of a Riemannian submersion.
\end{proof}

\section{Symplectic leaves}\label{symp-leaves}

It is well known that any smooth Poisson manifold foliates into smooth
symplectic initial submanifolds. These symplectic submanifolds are the
leaves of the integrable distribution with jumping rank defined by the
Hamiltonian vector fields. In particular, if $(M,\om)$ is a symplectic
\mf which is acted upon in a Hamiltonian fashion by a compact Lie
group $K$ such that this action is free this is true for the quotient
manifold $M/K$ which carries an induced Poisson structure. Moreover,
the symplectic leaves of $M/K$ are simply the connected components of
the Marsden-Weinstein reduced spaces $M\sporb K$ where $\orb$ is a
coadjoint orbit. 

In the present context where $M/K=(T^*Q)/K$
is a singular Poisson space  
such a statement cannot
be true in this form. 
(We continue to assume that $Q$ is of single isotropy type,
i.e., $Q=Q_{(H)}$ whence the cotangent lifted action has a non-trivial
isotropy lattice if $H\neq\set{e}$.) 
However, using that Hamiltonian vector fields
associated  to $K$-invariant functions are tangent to isotropy type
submanifolds $(T^*Q)_{(L)}$ we can still define a characteristic
distribution on $(T^*Q)_{(L)}/K$, and find in Subsection \ref{sub:sl}
that its symplectic leaves
are the connected components of the smooth symplectic manifold
$T^*(Q/K)\times_{Q/K}(\bsc_{q\in Q}\orb\cap\ann\ko_q)_{(L)}/K$. 
This result is put into perspective in Subsection \ref{sub:sp_reduc}
where we recall the singular Weinstein space description of
$T^*Q\sporb K \cong \WW\sporb K$.

\subsection{Universal reduction procedure}\label{sub:urp}
The singular reduction diagram of Ortega and Ratiu \cite[Theorem 8.4.4]{OR04}
adjoined to the universal reduction procedure of Arms, Cushman, and
Gotay \cite{ACG91}, see also \cite[Section 10.3.2]{OR04} applied to
the Weinstein space has the following form. 
\[
\xymatrix{
 {\mu_{A}^{-1}(\orb)}
   \ar @{<-^{)}}[r]
   \ar @{->>}[d]&
 {\mu_{A}^{-1}(\lam)}
   \ar @{^{(}->}[r]
   \ar @{->>}[d]&
 {\WW}
   \ar @{->>}[d]\\
 {\mu_{A}^{-1}(\orb)/K}
   \ar @{<-}[r]^-{\simeq}&   
 {\mu_{A}^{-1}(\lam)}/K_{\lam}
   \ar @{^{(}->}[r]&
 {\WW/K}
   \ar @{=}[r]&
 W  
}
\]
where $\lam\in\mu_{A}(\WW)$ and $\orb$ is the coadjoint orbit passing
through $\lam$. 
In this diagram the \iso $\mu_{A}^{-1}(\lam)/K_{\lam}\cong\WW\sporb K$
is an \iso of singular symplectic spaces with smooth structure, and
the inclusion  $\mu_{A}^{-1}(\lam)/K_{\lam}\hookto W$ is a morphism of
singular Poisson spaces. In the universal reduction scheme of Arms, Cushman, and
Gotay \cite{ACG91} this morphism is actually used to equip
$\mu_{A}^{-1}(\lam)/K_{\lam}$ with a Poisson structure.
Therefore, it is sensible 
to expect the smooth symplectic leaves of $W$ to be the connected
components of the smooth symplectic strata of $\WW\sporb K$. The
latter space is described in the next subsection.

\subsection{Symplectic reduction of $\WW$}\label{sub:sp_reduc}
Let $\orb$ be a coadjoint orbit in the image of the momentum
map $\mu_{A}: \WW\to\ko^{*}$,
and let $(L)$ be in the isotropy lattice of the $K$-action on
$\WW$ such that
$\lorb{\WW}:=\mu_{A}^{-1}(\orb)\cap\WW_{(L)}\neq\emptyset$. Then we have 
\[
 \lorb{\iota}: \lorb{\WW}\hookto\WW,
\]
the canonical embedding, and the orbit projection mapping
\[
 \lorb{\pi}: \lorb{\WW}\toto\lorb{\WW}/K=:(\WW\sporb K)_{(L)}.
\]
Consider furthermore 
\[
 \lorb{\rho}:
 (\bscorb)_{(L)}\toto(\bscorb)_{(L)}/K
\]
and 
\[
 \lorb{\nu}: 
 (\bscorb)_{(L)}\longto\orb,
 \quad
 (q,\lam)\longmapsto\lam
\]
as well as the embedding
\[
 \lorb{j}: 
 (\bscorb)_{(L)}\hookto\bsc_{q\in Q}\ann\ko_{q}.
\]
Finally, we denote the Kirillov-Kostant-Souriou symplectic form on
$\orb$ by $\omkks$, that is
$\omkks(\lam)(\ad^{*}(X).\lam,\ad^{*}(Y).\lam) = \vv<\lam,[X,Y]>$. 
Remember from Proposition \ref{prop:WW} that the symplectic structure
on $\WW$ is denoted by $\sigma$.

\begin{theorem}[Gauged symplectic reduction]\label{thm:WWorb}
Let $Q=Q_{(H)}$, 
let $\orb$ be a
coadjoint orbit in the image of the \momap $\mu_{A}: \WW\to\ko^{*}$,
and let $(L)$ be in the isotropy lattice of the $K$-action on
$\WW$ such that
$\lorb{\WW}:=\mu_{A}^{-1}(\orb)\cap\WW_{(L)}\neq\emptyset$. Then the
following are true.
\begin{enumerate}[\up (i)]
\item 
The smooth manifolds $(\WW\sporb K)_{(L)}$ and
  \[
   (\orb\spr{0}H)_{(L_{0})^{H}} =: (\orb\cap\ann\ho)_{(L_{0})^{H}}/H
  \]
  are typical symplectic strata of the stratified symplectic spaces
  $\WW\sporb K$ and $\orb\spr{0}H$ respectively. Here $L_{0}$ is an
  isotropy subgroup of the induced $H$-action on $\orb$ and
  $(L_{0})^{H}$ denotes its isotropy class in $H$. 
\item 
The symplectic stratum $(\WW\sporb K)_{(L)}$ can be globally
  described as
\[ 
  (\WW\sporb K)_{(L)}
  =
  T^{*}(Q/K)\times_{Q/K}(\bscorb)_{(L)}/K
\]
whence it is the total space of the smooth symplectic fiber bundle
\[
\xymatrix{
 {(\orb\spr{0}H)_{(L_{0})^{H}}}
 \ar @{^{(}->}[r]&
 {(\WW\sporb K)_{(L)}}
 \ar[r]
 &{T^{*}(Q/K)}
}
\]
  Hereby $L_{0}$ is an
  isotropy subgroup of the induced $H$-action on $\orb$ which is
  conjugate in $K$ to $L$, and
  $(L_{0})^{H}$ denotes its isotropy class in $H$. 
\item 
The symplectic structure $\lorb{\sigma}$ on $(\WW\sporb K)_{(L)}$
  is uniquely determined and given by the formula
\[
 (\lorb{\pi})^{*}\lorb{\sigma}
 =
 (\lorb{\iota})^{*}\sigma - (\mu_{A}|\lorb{\WW})^{*}\omkks.
\]
More precisely, 
\[
 \lorb{\sigma}=\Om^{Q/K}-\lorb{\beta}
\]
where $\lorb{\beta}\in\Om^{2}((\bscorb)_{(L)}/K)$ is defined by 
\[
 (\lorb{\rho})^{*}\lorb{\beta}
 =
  (\lorb{j})^{*}dB + (\lorb{\nu})^{*}\omkks.
\]
Finally $B$ is the form that was introduced in Proposition
\ref{prop:WW}. Thus for $(q,\lam)\in(\bscorb)_{(L)}$ and
\[
 (v_{i},\ad^{*}(X_{i}).\lam)\in T_{(q,\lam)}(\bscorb)_{(L)}
\]
where $i=1,2$ we have the
explicit formulas 
\[
 B_{(q,\lam)}(v_{i},\ad^{*}(X_{i}).\lam)
 = 
 \vv<\lam,A_{q}(v_{i})>
\]
and also
\begin{align*}
 d&B_{(q,\lam)}((v_{1},\ad^{*}(X_{1}).\lam),(v_{2},\ad^{*}(X_{2}).\lam))\\
 &=
 \vv<\lam,\curv_{q}^{A}(v_{1},v_{2})>
 + \vv<\lam,[X_{2},Z_{1}]>
 - \vv<\lam,[X_{1},Z_{2}]> 
 + \vv<\lam,[Z_{1},Z_{2}]>
\end{align*}
where $v_{i} = 
  \zeta_{Z_{i}}(q)\oplus v_{i}^{\textup{hor}}
  \in
  \ver_{q}\oplus\hor_{q}$ 
is the decomposition into vertical and horizontal parts with
  $Z_{i}\in\ko$.
\item  
The stratified symplectic space can be globally described as 
\[ 
 \WW\sporb K
 = 
 T^{*}(Q/K)\times_{Q/K}\bscorb/K
\]
whence it is canonically the total space of
\[
\xymatrix{
 {\orb\spr{0}H}
 \ar @{^{(}->}[r]&
 {\WW\sporb K}
 \ar[r]
 &{T^{*}(Q/K)}
}
\]
which is a singular symplectic fiber bundle with singularities
confined to the fiber direction in the sense of Definition \ref{def:sg-sp-fb}.
\end{enumerate}
\end{theorem}

\begin{proof}
Assertion (i) is a well-known  principle of stratified
symplectic reduction -- see Ortega and Ratiu \cite[Section 8.4]{OR04}.
The other assertions are proved in \cite[Theorem 5.5]{H04}
\end{proof}

\subsection{Symplectic leaves of $W=\WW/K$}\label{sub:sl}
Let $\orb$ be a coadjoint orbit in the image of the cotangent bundle
\momap $\mu: T^*Q\to\ko^*$. As before $(L)$ denotes an isotropy type
of the $K$-action on $T^*Q$, and $(L_0)^H$ denotes an isotropy type of
the $\Ad^*(H)$-action on $\orb$. 
Let $\lam\in\orb\cap\ann\ho$ with $K_{\lam}\cap H = L_0\subeq H$. We
want to make use of the Witt-Artin decomposition and thus denote the
symplectic slice of the $H$-action on $\orb$ at $\lam$ by 
\[
 V 
 :=
 (\ho.\lam)^{\Om^{\mathcal{O}}}/(\ho.\lam)
\] 
where $\Om^{\mathcal{O}}$
denotes the positive KKS-form on $\orb$. Notice that this is well
defined since $H.\lam=\Ad^*(H).\lam\subeq\orb\cap\ann\ho$ is an
isotropic submanifold of $\orb$. (The $H$-action on $\orb$ is
Hamiltonian with \momap $\orb\to\ho^*$ given by restriction to $\ho$.)
By construction $V$ is a symplectic vector space. Thus this is also
true for the linear subspace 
$V_{L_{0}} := \set{v\in V: H_v=L_0}$ 
of fixed symmetry type.   

\begin{lemma}\label{lem:V}
Under these assumptions 
$(\orb\cap\ann\ho)_{(L_0)^H}$ 
is a smooth manifold and
$
 T_{\lam}(\orb\cap\ann\ho)_{(L_0)^H} 
 = 
 T_{([e],0)}(H/L_0\times V_{L_{0}})
$.
Furthermore,
$(\orb\cap\ann\ho)_{(L_0)^H}/H$ 
is a smooth symplectic manifold and
$
 T_{[\lam]}((\orb\cap\ann\ho)_{(L_0)^H}/H) 
 = 
 V_{L_{0}}
$.
\end{lemma}

\begin{proof}
This is a direct consequence of the Symplectic Slice theorem. See Ortega and
Ratiu \cite[Chapter 7]{OR04} for a treatment of this theorem and
\cite[Section 8.1]{OR04} for the way in which it is used.
\end{proof}

For notational convenience we  abbreviate
$(\orb\cap\ann\ho)_{(L_0)^H}/H
 =: (\orb\spr{0}H)_{(L_0)^H}$.

\begin{lemma}\label{lem:span}
Assume that $[(q,\eta,\lam)]\in W_{(L)} = \WW_{(L)}/K$ and
$f\in\cinf(W)$ with unique lift $F=\phi^*f$.
Let $q\in Q_{H}$ and $\lam\in(\ann\ho)_{(L_0)^H}$ where $L_0\subeq H$
is conjugate to $L$ within $K$. Then 
\[
 \nabla^W_f[(q,\eta,\lam)]
 \in
 T_{([q],\eta)}(T^*(Q/K))\times T_{[\lam]}(\orb\spr{0}H)_{(L_0)^H}.
\]
Moreover, the latter space is the real span of local Hamiltonian
vector fields evaluated at $[(q,\eta,\lam)]$. 
\end{lemma}

\begin{proof}
According to Corollary \ref{cor:hvf}
we have 
\[
 \nabla_f^W[(q,\eta,\lam)]
 =
 (v_0,
  \left.(\Omega_{([q],\eta)}^{Q/K} \check{)}
         \right.^{-1}(d_{\widetilde{A}}F(q,\eta,\lambda)),
  \ad^*(d_{v}F(q,\eta,\lam)).\lam)
\]
where we consider $d_{\widetilde{A}}F(q,\eta,\lambda)$ as an element
of $T^*_{([q],\eta)}(T^*(Q/K))$ through the isomorphism given by the
dual of the horizontal lift with respect to the mechanical connection
$A$ on $Q\toto Q/K$. Moreover, 
$v_0 := 
 T_{([q],\eta)}\tau.\left.(\Omega_{([q],\eta)}^{Q/K} \check{)}
         \right.^{-1}(d_{\widetilde{A}}F(q,\eta,\lambda))$.
By Lemma \ref{lem:V} we have to check that $\ad^*(-Z).\lam\in V_{L_0}
= ((\ho.\lam)^{\Om^{\mathcal{O}}}/(\ho.\lam))_{L_0}$. Indeed, 
by Lemma \ref{lem3} it is 
true that $\Ad(h).Z=Z$ for all $h\in H$ where
$Z := -d_{v}F(q,\eta,\lam)$. Thus $[Z,\ho]=0$ whence
$\ad^*(-Z).\lam\in(\ho.\lam)^{\Om^{\mathcal{O}}}$ and even 
$\ad^*(-Z).\lam\in(\ho.\lam)^{\Om^{\mathcal{O}}}/(\ho.\lam)$. 
Moreover, it follows that 
$h.\ad^*(-Z).\lam = \ad^*(-Z).\lam$ for all $h\in H\cap
K_{\lam}=L_0$. However, since $\ad^*(-Z).\lam$ is also an element of
the Riemannian slice of the $H$-action on $\orb$ at $\lam$ it follows
that $H_{\ad^*(-Z).\lam}\subeq H_{\lam} = L_0$ since this is a general
feature of Riemannian slices. Therefore, $\ad^*(-Z).\lam\in V_{L_0}$,
as claimed.

The second claim follows by going again through the proof of Lemma
\ref{lem3}.
\end{proof}

Let $P_{(L)}$ denote the coinduced Poisson two-tensor on
$W_{(L)}$. Then we may rephrase the above lemma by saying that
\[
 \check{P}_{(L)}(T_{[(q,\eta,\lam)]}^*W_{(L)})
 =
 T_{([q],\eta)}(T^*(Q/K))\times T_{[\lam]}(\orb\spr{0}H)_{(L_0)^H}.
\] 
Thus we get the following theorem.

\begin{theorem}[Symplectic leaves]\label{thm:sp-leaf}
Let $(L)$ be an element of the isotropy lattice of the $K$-action on
$\WW$. Then the characteristic distribution of the coinduced Poisson
structure
$P_{(L)}$ on $W_{(L)}$ is given by
\[
 \check{P}_{(L)}(T^*W_{(L)})
 =
 T(T^*(Q/K)\times_{Q/K}(\bsc_{q\in Q}\orb\cap\ann\ko_q)_{(L)}/K)
\] 
whence the smooth symplectic leaves of $W_{(L)}$ are the connected
components of the smooth symplectic manifolds
\[
 (\WW\sporb K)_{(L)} 
 =
 T^*(Q/K)\times_{Q/K}(\bsc_{q\in Q}\orb\cap\ann\ko_q)_{(L)}/K.
\]
The symplectic form which makes the inclusion 
$(\WW\sporb K)_{(L)}\hookto W_{(L)}$ a Poisson morphism as that of 
Theorem \ref{thm:WWorb}.  
\end{theorem}

\begin{proof}
The statement about the integrability of the characteristic
distribution is tautologous since it is described as the tangent bundle
of a smooth manifold. The inclusion 
$(\WW\sporb K)_{(L)}\hookto W_{(L)}$ is a Poisson morphism by the
reasoning of Subsection \ref{sub:urp}.
\end{proof}

\subsection{Charge and spin}
We shortly describe and interpret some special cases of Theorem
\ref{thm:sp-leaf}. Retaining assumptions and notation of this theorem
we additionally require that $(\bsc_{q\in Q}\ann\ko_{q})/K\to Q/K$ be a trivial vector bundle so
that
\[
 (\bsc_{q\in Q}\orb\cap\ann\ko_q)_{(L)}/K)
 =
 Q/K\times \orb\spr{0}H
 \text{ and }
 \WW\sporb K  
 = 
 T^*(Q/K)\times\orb\spr{0}H.
\]
If $\orb\spr{0}H = \set{\textup{point}} = \set{\lam}$ then we obtain
the same result as in Corollary \ref{cor:charge}, i.e.\ $T^*(Q/K)$
inherits the magnetic symplectic form 
$\Om^{\lam} :=
\Om^{Q/K}-\vv<\lam,\tau^*\curv^{A}_{0}>$. 

If $\orb\spr{0}H$ is non-trivial then $\WW\sporb K  
 = 
 T^*(Q/K)\times\orb\spr{0}H$ with its induced symplectic structure
\[ 
 \Om(q,\eta,[\lam]) 
 =
 \Om^{Q/K}_{(q,\eta)}-\vv<\lam,\tau^*\curv^{A}_{0}>-\Om_{[\lam]}^{\mathcal{O}}
\]
(where $\Om^{\mathcal{O}}$ denotes the reduced symplectic form on
$\orb\spr{0}H$ obtained from the KKS-form on $\orb$)  describes the
phase space of a (color)-charged  particle moving on $Q/K$ under the
influence of the field $\curv^A_0$ and with internal spin parameters
corresponding to $\orb\spr{0}H$.

\section{Examples}\label{sec:ex}

\subsection{Calogero-Moser space}
In the spirit of 
Hochgerner \cite{H04} 
we can apply Theorem \ref{thm:p-struct} to obtain rational versions of
spin Calogero-Moser systems. That is,
let $V$ be a real Euclidean
vector space and $K$ a connected compact Lie group that acts on $V$ through
a polar representation. Via the inner product we consider the
cotangent bundle of $V$ as a product $T^{*}V=V\times V$. The canonical
symplectic form $\Om$ is thus 
\[
 \Om_{(a,\alpha)}((a_{1},\alpha_{1}),(a_{2},\alpha_{2})) =
 \vv<\alpha_{2},a_{1}>-\vv<\alpha_{1},a_{2}>
\]
where $\vv<\phantom{a},\phantom{a}>$ is the inner product on $V$. 

The
cotangent lifted action of $K$ is the diagonal action
on $V\times V$. According to Dadok \cite{Dad85} we may 
think of the action by $K$ on $V$ as a symmetric space representation
and thus consider
$\ko\oplus V =: \lo$ as a real semisimple Lie algebra with Cartan decomposition
into $\ko$ and $V$, and with bracket relations
$[\ko,\ko]\subeq\ko$, $[\ko,V]\subeq V$, and $[V,V]\subeq V$. The
momentum mapping corresponding to the $K$-action on $T^{*}V=V\times V$
is now given by 
$
 \mu: V\times V\to\ko^{*}=\ko$,
$
 (a,\alpha)\mapsto[a,\alpha]=\ad(a).\alpha
$
where we identify $\ko=\ko^{*}$ via an $\Ad(K)$-invariant inner product. 

Let $V_r$ denote the open and dense subset of $V$ consisting of
regular elements.
Following the ideas of \aklm \cite{AKLM03} we can interpret the
Poisson reduced space $(T^*V_r)/K$ as the Poisson phase space of a
rational spin Calogero-Moser system. Indeed, let $H$ denote the free
Hamiltonian on $T^*V_r$ and $h$ its induced Hamiltonian on
$(T^*V_r)/K$. Then the Hamiltonian \vf of $h$ computed by Theorem
\ref{thm:p-struct} and its Corollary \ref{cor:hvf}
is exactly minus the one of the Calogero-Moser Hamiltonian
in \cite[Section 6.5]{AKLM03}. (The difference in the sign comes from
a different convention in defining the Hamiltonian \vf of a function.)
According to Section \ref{symp-leaves} the smooth symplectic leaves of
$(T^*V_r)/K$ are given by the connected components of smooth strata of
$(T^*V_r)\sporb K$ where $\orb$ is a coadjoint orbit in $\ko^*$. The
latter space is described in \cite[Section 6.3]{H04}.

\subsection{Orbit type $SO(5)/SO(3)$}
Let $SO(5)$ act on $S^9\subeq\R^5\times\R^5$ through the diagonal
action. We denote by $Q$ the open and dense subset of elements $(v,w)\in
S^9\subeq\R^5\times\R^5$ such that $v$ and $w$ are linearly
independent. Clearly, $Q$ is preserved by the $SO(5)$-action, and constitutes,
moreover, the regular stratum with respect to this action. Thus $Q$ is
of single isotropy type, and this type is easily seen to be
$H := SO(3)\subeq SO(5) =:K $. Writing $K$ as a matrix group we embed
$H$ in the usual way in the lower right corner. The orbit space $Q/K$
can be diffeomorphically identified with the open disk $B^2$ of radius
$1$ in  $\R^2$, and the
projection $Q\toto Q/K$ is a (non-principal) fiber bundle with typical
fiber $K/H$. 

We consider the cotangent lifted $K$-action on $T^*Q$. 
This is clearly a non-free action with a non-trivial isotropy lattice.
By Theorem \ref{thm:p-struct} the singular Poisson reduced space with
respect to the lifted $K$-action is of the form
\[
 T^*B^2\times\ann\ho/\Ad^*(H),
\]
since the bundle $(\bsc_{q\in Q}\ann\ko_q)/K\to Q/K=B^2$ is, in this case, trivial. 
Using the trace
form we identify $\ann\ho$ with $\ho^{\bot}$. Now, the map 
\[
 \ho^{\bot}\longto\R\times\R^3\times\R^3,
 \text{ }
 (x_{ij})_{ij}\longmapsto(x_{21},(x_{k1})_{k=3}^{5},(x_{k2})_{k=3}^{5})
  = (t,v,w)
\]
is a linear isomorphism that is equivariant for the $H$-action on the
right hand side which acts trivially on the $\R$-factor and in the
standard diagonal way on the $\R^3\times\R^3$-factor.
Thus the singular Poisson reduced space with
respect to the lifted $K$-action is of the form
\[
 T^*B^2\times\R\times\R^3\times_H\R^3.
\]
However, the induced Poisson structure is not obvious at all
(if we did not have Theorem \ref{thm:p-struct}). 
The stratification is the product stratification induced by the obvious
stratification of $(\R^3\times\R^3)/H$. 
Notice also that the induced form $\curv_0^A$ on $B^2$ which comes
from the mechanical curvature is by Proposition \ref{prop:curvA}
$\R$-valued and evaluates on the $t$-factor. 

In order to get a non-trivial symplectic leaf of
$T^*B^2\times\ann\ho/\Ad^*(H)$ let 
\[
 \lam
 := 
\left(
 \begin{matrix}
   0 & 1 & 1 & 0 & 0\\
  -1 & 0 &-1 & 1 & 0\\
  -1 & 1\\
   0 &-1\\
   0 & 0\\
 \end{matrix}
\right)
\]
and $\orb$ be the (co-)adjoint orbit through $\lam$. Doing the
appropriate linear algebra one sees that $\dim\ko_{\lam} = 2$,
$\dim\ko_{\lam}^{\bot}=8$, $\ho\cap\ko_{\lam}=\set{0}$, and
$\dim\ho^{\bot}\cap\ko_{\lam}^{\bot} = 5$. Doing a little more linear
algebra the symplectic normal space (see, in particular, 
Subsection \ref{sub:sl})
to the $H$-action on $\orb$ at $\lam$ computes to be 
\[
 V =
\set{
\left(
 \begin{matrix}
   0 & t & a & 0 & 0\\
  -t & 0 & t & a & 0\\
  -a &-t\\
   0 &-a\\
   0 & 0\\
 \end{matrix}
\right)
 : t,a\in\R}.
\]
Since $\ho\cap\ko_{\lam}=\set{0}$ we thus get $V_{L_0} = V$ whence 
$
 T_{[\lam]}((\orb\cap\ann\ho)_{(L_0)^H}/H) 
 = 
 V
$.
In particular, the symplectic leaf passing through 
\[ 
 T^*B^2\times\set{[\lam]}
  \subeq T^*B^2\times\ann\ho/\Ad^*(H) 
  =  T^*B^2\times\R\times\R^3\times_H\R^3
\]
is $6$-dimensional. Further and more detailed 
investigation into this example is
written up in \cite{h05}. 
It actually turns out that the mechanical curvature yields a
magnetically non-trivial 
symplectic structure on $T^*B^2$.

\end{document}